\documentclass[10pt,a4paper]{article}
\usepackage[margin=2.5cm]{geometry}\usepackage[utf8]{inputenc}
\usepackage[T1]{fontenc}
\usepackage{lmodern}
\usepackage{amsmath,amsfonts,amsthm,amssymb,mathtools}
\usepackage{xcolor,microtype,graphicx,caption}
\usepackage{tikz}
\usetikzlibrary{arrows.meta,decorations.markings}
\usepackage{hyperref}
\colorlet{darkblue}{blue!50!black}
\hypersetup{
  colorlinks=true,
  linkcolor=darkblue,
  citecolor=blue,
  urlcolor=blue,
  filecolor=red,
  pdfstartview={FitH}
}
\theoremstyle{plain}
\newtheorem{theorem}{Theorem}[section]
\newtheorem{lemma}[theorem]{Lemma}
\newtheorem{proposition}[theorem]{Proposition}
\newtheorem{corollary}[theorem]{Corollary}
\newtheorem{theoremA}{Theorem}

\newtheorem{theoremB}{Theorem}

\theoremstyle{definition}
\newtheorem{definition}[theorem]{Definition}
\theoremstyle{remark}
\newtheorem{remark}[theorem]{Remark}
\theoremstyle{remark}

\theoremstyle{plain}

\newcommand{\R}{\mathbb{R}}
\newcommand{\C}{\mathbb{C}}
\newcommand{\Z}{\mathbb{Z}}

\newcommand{\T}{\mathbb{T}}
\newcommand{\HH}{\mathcal{H}}
\newcommand{\cR}{\mathcal{R}}

\newcommand{\bS}{\mathbb{S}}

\newcommand{\cM}{\mathcal{M}}
\newcommand{\cF}{\mathcal{F}}
\newcommand{\cJ}{\mathcal{J}}

\newcommand{\cN}{\mathcal{N}}
\newcommand{\cZ}{\mathcal{Z}}

\newcommand{\cO}{\mathcal{O}}
\newcommand{\Id}{\mathrm{Id}}
\newcommand{\dd}{\mathrm{d}}
\newcommand{\e}{\varepsilon}
\newcommand{\lan}{\langle}
\newcommand{\ran}{\rangle}
\DeclareMathOperator{\Real}{Re}
\DeclareMathOperator{\Imag}{Im}

\DeclareMathOperator{\sgn}{sgn}
\numberwithin{equation}{section}
\numberwithin{figure}{section}

\title{Small-time global exact controllability
of nonlinear Schr\"odinger equations on the circle
with bilinear controls}
\author{Vahagn Nersesyan\thanks{NYU-ECNU Institute of Mathematical Sciences at NYU Shanghai, 3663 Zhongshan Road North, Shanghai, 200062, China; e-mail: vahagn.nersesyan@nyu.edu}}
\date{\today}
\begin{document}
\maketitle

\begin{abstract}
We study nonlinear Schr\"odinger equations on the circle with
nonlinearity $\kappa|\psi|^{2p}\psi$, a fixed smooth real potential
$V$, and bilinear controls. Using only the three control potentials
$\mathbf 1,\cos x,\sin x$, we prove global approximate
controllability in $H^s$ on the unit sphere of $L^2$ for all
integers $s,p\ge1$, every $\kappa\in\mathbb R$, and any prescribed
positive time. The proof combines two ingredients: the approximate
controllability between nowhere vanishing states from our companion
paper, and a new mechanism that steers any nonzero initial state to
a nowhere vanishing state in arbitrarily small time.

When $V=0$ and the nonlinearity is defocusing or sufficiently weakly
focusing, adding two explicit control potentials of finite
regularity gives local exact controllability in $H^2$ near the
constant state, via a moment problem and the inverse mapping
theorem. Together, these results imply global exact controllability
in $H^2$ on the unit sphere of $L^2$, at any prescribed positive
time, with five control potentials. To~the best of our knowledge,
this is the first global exact controllability result for a
nonlinear Schr\"odinger equation with finitely many bilinear
controls.
\end{abstract}

\medskip\noindent
\textbf{2020 Mathematics Subject Classification:}
Primary 93B05; Secondary 35Q55, 93C20, 81Q93.

\noindent
\textbf{Keywords:} Nonlinear Schr\"odinger equation, bilinear control,
polynomial nonlinearity, global exact controllability,
removal of zeros, moment problem.

\setcounter{tocdepth}{2}
\tableofcontents
\setcounter{tocdepth}{2}

\section*{Introduction}
\addcontentsline{toc}{section}{Introduction}

In this paper, we study the controllability of the nonlinear
Schr\"odinger equation on the circle $\T=\R/2\pi\Z$
\begin{equation}\label{E:0.1}
i\partial_t\psi=-\partial_x^2\psi+V(x)\psi
+\kappa|\psi|^{2p}\psi+\lan u(t),Q(x)\ran\psi.
\end{equation}
Here $p\ge1$ is an integer, $\kappa\in\R$, and $V$ is a smooth
real-valued potential. The control $u\in L^2(0,T;\R^q)$ acts
through the time-dependent potential
$$
\lan u(t),Q(x)\ran=\sum_{j=1}^q u_j(t)Q_j(x),
$$
where $Q_1,\dots,Q_q$ are fixed real-valued functions.
We supplement equation~\eqref{E:0.1} with the initial condition
\begin{equation}\label{E:0.2}
\psi(0,x)=\psi_0(x),
\end{equation}
where $\psi_0\in H^s=H^s(\T;\C)$ and $s\ge1$ is an integer.
For smooth control potentials, the Cauchy
problem~\eqref{E:0.1},~\eqref{E:0.2} is locally well posed in $H^s$.

Since the potentials are real-valued, the $L^2$ norm of the
solution is conserved.  We therefore consider
initial and target states on the unit sphere
$$
\bS=\{\psi\in L^2(\T;\C):\|\psi\|_{L^2}=1\}.
$$
Our main result is global exact controllability in $H^2$ at any
prescribed positive time. More precisely, when~$V=0$ and the
nonlinearity is defocusing or sufficiently weakly focusing, we show
that any initial state in $H^2\cap\bS$ can be steered exactly to any
target state in $H^2\cap\bS$ by using five explicit control
potentials. Three of them are smooth and two are of finite
regularity. These potentials do not depend on the initial and target
states, on the time $T>0$, or on the admissible parameters $p$ and
$\kappa$. Only the controls $u_j(t)$, $j=1,\ldots,5$, depend on these
data. For example, the constant state $\phi_0=(2\pi)^{-1/2}$ can be
steered to any state in $H^2\cap\bS$ supported in an arbitrarily
short arc. Thus, in spite of the nonlinearity and dispersion of the
equation, the entire mass can be concentrated in this arc in
arbitrarily small time. When the target state is not smooth, this is
impossible with smooth control potentials, since a solution starting from~$\phi_0$ remains smooth as long as it
exists. This is the reason why we
use control potentials of finite regularity.

We first prove a global approximate controllability result in
$H^s$, which is of independent interest and holds under more
general assumptions. It requires only the three potentials
$\mathbf1,\cos x,\sin x$ and holds for all integers $s,p\ge1$,
every smooth real-valued potential $V$, and every $\kappa\in\R$,
with no restriction on the sign or the size of the nonlinearity.
Let us now formulate these results in more detail.

\subsection*{Global approximate controllability}

Let $Q_1,\dots,Q_q$ be smooth real-valued functions such that
\begin{equation}\label{E:0.3}
\{\mathbf1,\cos x,\sin x\}
\subset\operatorname{span}\{Q_j:j=1,\dots,q\}.
\end{equation}
This condition guarantees that the controls generate a sufficiently
rich family of phases, as it implies the saturation property
\textup{(\hyperlink{H1}{H$_1$})} of Section~\ref{S:1}. Our global approximate
controllability result reads~as~follows.

\begin{theoremA}\label{T:A}
Assume~\eqref{E:0.3} and let $s\ge1$ be an integer.
For any $\psi_0,\psi_1\in H^s\cap\bS$, any $T>0$, and any~$\e>0$, there are a control $u\in L^2(0,T;\R^q)$ and a solution
$\psi\in C([0,T];H^s\cap\bS)$ of~\eqref{E:0.1} such that
$$
\psi(0)=\psi_0,\qquad \|\psi(T)-\psi_1\|_s<\e.
$$
\end{theoremA}

As mentioned above, this theorem holds for every smooth real-valued
potential $V$, every integer $p\ge1$, and every $\kappa\in\R$. Note
that the control potentials $\mathbf 1,\cos x,\sin x$ do not depend
on these parameters, and that the time $T$ can be arbitrarily small.
Even when solutions of the uncontrolled equation blow up in finite
time, the controlled solution given by the theorem exists on the
whole interval $[0,T]$.

In the companion paper~\cite{NLSI}, we establish approximate 
controllability in $H^s$ of the nonlinear
Schr\"odinger equation  on a torus $\T^d$ of arbitrary dimension,
for every integer $s>d/2$, between states that do not vanish. On the circle, such states are dense, and
the restriction on the target can be removed by approximation. The
restriction on the initial state is much more serious. The controls act
by multiplication: they change the phase of the state, but not its
modulus. Thus they cannot, by themselves, create mass at the points
where the state vanishes. The main novelty of
Theorem~\ref{T:A} is a mechanism that removes the zeros of any
nonzero initial state in arbitrarily small time.

This mechanism is based on the dispersion of the free Schr\"odinger
equation. At suitable short times, the free flow on the circle
transforms the state into a superposition of finitely many shifted
copies of itself. Thus every point of the circle receives
contributions from several parts of the initial state, including
parts where it does not vanish. Note that these contributions may 
cancel each other. To avoid this, we first multiply the initial
state by a suitable phase. As shown in~\cite{DN}, this phase multiplication can be approximated in
$H^s$ by solutions of the controlled nonlinear equation. We prove
that the same is true for the short free evolution, for any smooth
potential $V$ and any $\kappa\in\R$. Since the time needed for these
two operations can be chosen arbitrarily small, we obtain a state
without zeros in arbitrarily small time. We then apply the result
of~\cite{NLSI} on the remaining time interval to approximate~the~target.

Let us note that this argument is specific to the circle, as is the
density of nowhere vanishing states. On~the~circle, the zeros of a state can be removed by an arbitrarily
small perturbation. In higher dimensions, this is not true in
general: some zeros persist under small perturbations;
see~\cite[Remark~3.3]{NLSI}. For~this reason, our argument does not extend
directly to the multidimensional case.

\subsection*{Global exact controllability}

Let us now assume that $V=0$. By Theorem~\ref{T:A}, any state
in $H^2\cap\bS$ can be steered arbitrarily close to the constant
state $\phi_0$ in $H^2$. Therefore, to prove global exact
controllability, it remains to establish local exact controllability
near $\phi_0$ in the same space. We then combine these two results, using the time
reversibility~of~the~equation.

As explained above, the local result cannot be obtained with smooth
control potentials, since every $H^2$ neighborhood of $\phi_0$
on~$\bS$ contains nonsmooth states. We therefore add to
$\mathbf1,\cos x,\sin x$ the two potentials
\begin{equation}\label{E:0.4}
Q_c(x)=\sum_{n=1}^\infty\frac{\cos nx}{n^2},
\qquad
Q_s(x)=\sum_{n=1}^\infty\frac{\sin nx}{n^2}.
\end{equation}
These potentials belong to $H^r(\T)$ for any $r<3/2$, but not to
$H^{3/2}(\T)$. They allow us to overcome the smoothness obstruction,
and a smoothing estimate for the Duhamel term still gives local
well-posedness in~$H^2$. To prove local exact controllability, we develop the approach
of~\cite{BL}, where a smoothing estimate for the Duhamel term is
combined with the inverse mapping theorem to prove local exact
controllability of Schr\"odinger equations with bilinear controls
on an interval, both for the linear equation with Dirichlet
boundary conditions and for the cubic equation with Neumann
boundary conditions.

Our result remains true for a class of perturbations of the
potentials $Q_c$ and $Q_s$. Let
$$
\mathcal Q=(Q_c,Q_s)+H^2(\T;\R)^2,
$$
endowed with the topology induced by translation from
$H^2(\T;\R)^2$. For
$(P_c,P_s)\in\mathcal Q$, we set
$Q=(\mathbf1,\cos x,\sin x,P_c,P_s)$ and consider the equation
\begin{equation}\label{E:0.5}
i\partial_t\psi=-\partial_x^2\psi+\kappa|\psi|^{2p}\psi
+\lan u(t),Q(x)\ran\psi.
\end{equation}
By a solution on $[0,T]$, we mean a function in $C([0,T];H^2)$
satisfying the Duhamel formula in $H^1$. With~this definition, the
problem is locally well posed in $H^2$.

\begin{theoremB}\label{T:B}
There is an open dense set $\cO\subset\mathcal Q$ containing
$(Q_c,Q_s)$ with the following property. Let~$(P_c,P_s)\in\cO$, and
let $p\ge1$ be an integer and $\kappa\in\R$ be such that
\begin{equation}\label{E:0.6}
\beta:=2p\kappa(2\pi)^{-p}>-4.
\end{equation}
Then, for any $T>0$ and any $\psi_0,\psi_1\in H^2\cap\bS$, there are
a control $u\in L^2(0,T;\R^5)$ and a solution
$\psi\in C([0,T];H^2\cap\bS)$ of~\eqref{E:0.5} such that
$$
\psi(0)=\psi_0,\qquad \psi(T)=\psi_1.
$$
\end{theoremB}

Let us note that the set $\cO$ does not depend on $T$, $p$, and
$\kappa$. Since $\cO$ is open and contains $(Q_c,Q_s)$, the theorem
applies, in particular, to all sufficiently small $H^2$ perturbations
of the explicit pair~\eqref{E:0.4}.

Condition~\eqref{E:0.6} is satisfied in the linear case $\kappa=0$,
for all defocusing nonlinearities $\kappa>0$, and for focusing
nonlinearities with sufficiently small $|\kappa|$. It is equivalent
to the inequality $n^2+\beta>0$ for all~$n\ge2$. Under this
condition, the equation linearized at $\phi_0$ oscillates on every
Fourier mode $n\ge2$, and the gaps between the corresponding
frequencies tend to infinity. This allows us to control these modes
in arbitrarily small time. The first mode is controlled directly by
$\cos x$ and $\sin x$, so no condition is needed for $n=1$. We do not
expect condition~\eqref{E:0.6} to be optimal. For $\beta\le-4$, some modes of the linearization at $\phi_0$ are no
longer oscillatory, and for certain values of $\beta$, distinct modes
resonate; see Remark~\ref{R:4.7}. This prevents us from using the
same linearized argument, but does not rule out nonlinear
controllability.

To the best of our knowledge, Theorem~\ref{T:B} is the first
global exact controllability result for a nonlinear Schr\"odinger
equation with finitely many bilinear controls. All the global exact
controllability results with bilinear controls that we are aware of
concern linear equations; see, for example,~\cite{Ner10,MN15}.

\subsection*{Previous work}

The controllability of Schr\"odinger equations with bilinear controls
has been intensively studied over the past two decades. One of the main
motivations comes from quantum physics and chemistry, where the
control is typically an external electric field acting on the system
through one or several fixed spatial profiles. Most of the literature
on this problem deals with the linear equation, that is, with the
case $\kappa=0$.

An important negative result is proved
in~\cite[Theorem~3.6]{BMS82}: for a linear evolution equation with
a bounded bilinear control operator on an infinite-dimensional
Banach space, the set of states reachable from a given initial
state with $L^2$ controls is contained in a countable union of
compact sets. In particular, since this set has empty interior,
such an equation is not exactly controllable in this space, not
even locally.

The first positive result on exact controllability of the
Schr\"odinger equation with bilinear controls was obtained
in~\cite{Bea05}, where local exact controllability near the ground
state of the linear equation on an interval with Dirichlet boundary
conditions is proved by the return method~\cite{Cor07} and a
Nash--Moser theorem. It was later shown in~\cite{BL} that a
smoothing estimate for the Duhamel term allows one to use the
classical inverse mapping theorem. Under some coupling assumptions,
this gives local exact controllability in any positive time for the
linear equation with Dirichlet boundary conditions and for the cubic
equation with Neumann boundary conditions. For a nonlinear equation
controlled through the length of the interval, local exact
controllability near the ground state is proved in~\cite{BLT15}.
More recently, local exact controllability of the cubic equation
with Dirichlet boundary conditions near the ground state was proved
in~\cite{DN-Dirichlet}, using four bilinear controls. For
the linear equation with Dirichlet boundary conditions, local exact
controllability in $H_0^1$ with discontinuous control potentials is
obtained in~\cite{BD}, where periodic $H^1$ states are also treated
by adding a constant magnetic potential that makes the spectrum
simple.

Our local result is closest to~\cite[Theorem~4]{BL}, where local
exact controllability of the cubic equation with Neumann boundary
conditions is proved in $H^2$ near a spatially constant solution. For $n\ge1$, the
assumption there is that the $n$th cosine coefficient of the control
potential has absolute value at least $c n^{-2}$. On an interval,
a smooth potential can satisfy this condition because integration by parts gives boundary terms. On the circle,
the Fourier coefficients of a smooth potential decay faster than
any inverse power of $n$. We therefore use the potentials $Q_c$
and $Q_s$ of finite regularity. Their cosine and sine coefficients,
respectively, are of order $n^{-2}$, and the two potentials allow
us to control the cosine and sine modes independently in each
double eigenspace. Multiplication by these potentials does not map~$H^2$ into itself. Therefore, the control operator is not bounded in
$H^2$, and the obstruction of~\cite{BMS82} does not apply in this
space. On the other hand, $Q_c$ and $Q_s$ are bounded functions, so,
in the linear case, the obstruction applies in $L^2$. This is one of
the reasons why we work in $H^2$. As in~\cite{BL}, we use a smoothing
estimate for the Duhamel term to construct solutions in $H^2$.
Unlike~\cite{BD}, we do not modify the equation: the nonzero
eigenvalues of the periodic Laplacian remain double, and we control
the corresponding cosine and sine modes with two different
potentials. A similar strategy was used in the parabolic setting
in~\cite{DPU25}.

Global approximate controllability has been studied by several
methods. A geometric approach based on finite-dimensional Galerkin
approximations was introduced in~\cite{CMSB09} and further developed
in many subsequent papers; see, for example,~\cite{BCCS12,BCS14} for the case of several
controls. In particular, it follows from these results that the
linear equation on the circle with~$V=0$ is globally approximately
controllable in~$L^2$ with the two potentials~$\cos x$ and~$\sin x$;
see~\cite[Proposition~8.2]{BCCS12}. Another approach is based on
Lyapunov techniques. In~\cite{Ner10}, such techniques are used to
prove approximate controllability to eigenstates in higher Sobolev
norms for a generic control potential. Combined with a local exact
result, this gives global exact controllability for regular states
on an interval. This result is extended to the simultaneous control
of several equations in~\cite{MN15}. All these results concern linear equations: they rely on the linear
structure of the equation and do not extend directly to the
nonlinear case. Moreover, these global results do not, in general,
give controllability at an arbitrary prescribed positive time.

More recent results are based on the Agrachev--Sarychev
method~\cite{AS05,AS06}, introduced for the Euler and Navier--Stokes
equations with additive controls acting on finitely many Fourier
modes. For the defocusing cubic Schr\"odinger equation on the
two-dimensional torus driven by an additive control of finite
dimension, approximate controllability and exact controllability of
the projections onto any finite number of Fourier modes are proved
in~\cite{Sar12}. In~\cite{DN}, the method is extended to bilinear
controls, and it is shown that, for any $\psi_0\in H^s$ and any
smooth real periodic function $\theta$, the state $e^{i\theta}\psi_0$
can be approximated in $H^s$ by controlled solutions starting from
$\psi_0$ in arbitrarily small time. Combining this phase control with
transport by diffeomorphisms, small-time global approximate
controllability in $L^2$, modulo a global phase factor, is obtained
in~\cite{BP} for linear equations and in~\cite{BCP} for the
logarithmic Schr\"odinger equation. Let us note that a diffeomorphism,
like a phase multiplication, does not remove the zeros of a state, but
only moves them.

Small-time approximate controllability has also been studied in
other geometries. On $\R^d$, space dilations and time contractions
are used in~\cite{BPexamples} to construct examples of linear
equations that are approximately controllable in small time, also in
stronger norms related to the harmonic oscillator. On compact
connected Riemannian manifolds without boundary, finite saturating
families of potentials are constructed in~\cite{BDPNash}, which gives
small-time global approximate controllability in $L^2$ for the linear
equation with $V=0$.

Finally, let us mention some results for additive controls localized
in a subdomain. Localization is another natural type of degeneracy of
the control, different from the finite dimensionality considered
here, and it is also important in applications. For the cubic Schr\"odinger equation on the circle, local exact
controllability near zero in $H^s$, for any $s\ge0$ and $T>0$, is
established in~\cite{RZ09}. For the same equation on an interval,
including the case of periodic boundary conditions, global exact
controllability in $L^2$ in sufficiently large time, uniformly on
bounded sets, is proved in~\cite{Lau10}, by combining stabilization
with local exact controllability near zero, as in~\cite{DGL06}. In our case, the bilinear
controls preserve the $L^2$ norm, so the local analysis is carried
out near the constant state $\phi_0$ instead of zero.

\subsection*{Outline of the proofs}

The main new step in the proof of global approximate controllability is the removal of zeros. As explained above, we multiply a nonzero initial state $\psi_0$ by a phase $e^{i\theta}$ and then let it evolve under the free Schr\"odinger equation for a short time $\tau>0$. Let us explain how $\theta$ and $\tau$ are chosen. We use the following classical property of the free Schr\"odinger group $U(t)=e^{it\partial_x^2}$ on the circle: for $\tau=2\pi/N$ with odd $N$, the operator~$U(\tau)$ is a finite sum of translations by $2\pi k/N$ with nonzero coefficients. If $N$ is large, then at every point of the circle at least three of the translated copies of $\psi_0$ are nonzero. The difficulty is to find a smooth phase $\theta$ that prevents the cancellation of these copies at all points at once. Using a Fourier argument, we show that, at a given point, the sum is small only for a small set of phases. To extend this to the whole circle, we divide it into short intervals and use the $H^1$ regularity of $\psi_0$ on each of them. Here the one-dimensional setting is essential: the number of intervals is small enough to be compensated by the pointwise estimate, which is no longer the case in higher dimensions. As a result, almost every phase in a suitable finite-dimensional family prevents the cancellation at all points.

It remains to reach a state close to $U(\tau)(e^{i\theta}\psi_0)$
with the controlled nonlinear equation. The multiplication by
$e^{i\theta}$ can be approximated by the phase control of~\cite{DN},
so the main point is to approximate the free evolution. On a short
time interval, the potential $V$ and the nonlinearity act essentially
as a phase. We therefore alternate short intervals of uncontrolled
evolution with phase controls that compensate for this action, in the
spirit of the Lie--Trotter product formula. In this way, the
controlled solution approximates the free evolution as closely as we
wish. Since $U(\tau)(e^{i\theta}\psi_0)$ has no zeros, the same is
true for any sufficiently close state, and we obtain a state without
zeros in arbitrarily small time. Starting from this state, we then
use the approximate controllability of~\cite{NLSI} for initial states
without zeros to steer the solution close to the target in the
remaining time.

Let us now turn to exact controllability, assuming that $V=0$. The
potentials $Q_c$ and $Q_s$ are singular only at $x=0$, and
multiplication by them does not preserve $H^2$. However, for $P=Q_c$
or $Q_s$ and~$v\in H^2$, the product $P(v-v(0))$ belongs to $H^2$,
and the Duhamel term corresponding to $v(0)P$ can be estimated in~$H^2$ by using the $n^{-2}$ decay of the Fourier coefficients of $P$.
This gives local well-posedness of the nonlinear equation in $H^2$
and smooth dependence of the solution on the initial state and the
control. The same decay is exactly what is needed to control the
equation linearized at the stationary state $\phi_0$ with targets in
$H^2$. Indeed, each nonzero eigenvalue of the periodic Laplacian
corresponds to a cosine mode and a sine mode, and, for $n\ge2$, we
control the cosine modes with~$Q_c$ and the sine modes with~$Q_s$.
If~$\beta>-4$, the frequencies of these modes are positive and their
gaps tend to infinity, so we can solve the corresponding moment
problem on any time interval. The first mode is controlled by
$\cos x$ and $\sin x$, and the constant phase by $\mathbf1$. As
in~\cite{BL}, we then apply the inverse mapping theorem and obtain
local exact controllability on the unit sphere near $\phi_0$ in any
positive time. Finally, by approximate controllability in $H^2$, we
steer the initial state to this neighborhood, and, using time
reversibility, we reach the target from this
neighborhood. The same proof works for an open dense set of $H^2$
perturbations of~$(Q_c,Q_s)$, for which the two potentials still act
independently at every frequency $n\ge2$.

The paper is organized as follows. In Section~\ref{S:1},
we recall some results on well-posedness and controllability that are
used throughout the paper. Section~\ref{S:2} is devoted to
global approximate controllability. We first show that the free
evolution can be approximated by solutions of the controlled
equation, and then we use this property to remove the zeros of the
initial state. In Section~\ref{S:3}, we establish
well-posedness in $H^2$ for the equation with the additional
potentials $Q_c$ and $Q_s$. In Section~\ref{S:4}, we prove local
exact controllability near $\phi_0$ and derive global exact
controllability. Finally, in Appendix~\ref{SS:A}, we
construct the phase used in the removal of zeros, and in
Appendix~\ref{SS:B}, we prove that the admissible pairs of
potentials form an open dense set.

\section{Preliminaries}\label{S:1}

In this section, we recall some results on local well-posedness and
controllability for equation~\eqref{E:0.1}. Let~$s\ge1$ be
an integer and assume that $V,Q_1,\dots,Q_q$ are smooth
real-valued functions.

For $k\ge0$, we endow the Sobolev space $H^k=H^k(\T;\C)$ with the norm
$$
\|f\|_k^2=2\pi\sum_{n\in\Z}(1+n^2)^k|\hat f(n)|^2,
\qquad
\hat f(n)=\frac1{2\pi}\int_\T f(x)e^{-inx}\,\dd x.
$$

Let $U(t)=e^{it\partial_x^2}$ denote the free Schr\"odinger group.
Given $T>0$, $\psi_0\in H^s$, and $u\in L^2(0,T;\R^q)$,
by a solution of equation~\eqref{E:0.1} with initial condition $\psi(0)=\psi_0$
and control $u$, we mean a function $\psi\in C([0,T];H^s)$
satisfying the Duhamel formula
\begin{equation}\label{E:1.1}
\psi(t)=U(t)\psi_0-i\int_0^t U(t-r)
\big[V\psi(r)+\kappa|\psi(r)|^{2p}\psi(r)
+\lan u(r),Q\ran\psi(r)\big]\,\dd r, \quad t\in[0,T].
\end{equation}
 The following local well-posedness result follows
from~\cite[Proposition~1.1]{DN}.

\begin{proposition}\label{P:1.1}
For every $\psi_0\in H^s$ and $u\in L^2_{\mathrm{loc}}(\R_+;\R^q)$,
there is a unique maximal solution $\psi\in C([0,T_*);H^s)$
of~\eqref{E:0.1} with $\psi(0)=\psi_0$, where $T_*\in(0,\infty]$.
If $T_*<\infty$, then $\|\psi(t)\|_s\to\infty$ as $t\uparrow T_*$.
For every $T\in(0,T_*)$, the solution map is defined and locally
Lipschitz on a neighborhood of $(\psi_0,u|_{(0,T)})$ in
$H^s\times L^2(0,T;\R^q)$, with values in $C([0,T];H^s)$.
Moreover,
$$
\|\psi(t)\|_{L^2}=\|\psi_0\|_{L^2},\qquad 0\le t<T_*.
$$
For every $R>0$, there is $\tau_R>0$ such that, for $u=0$
and $\|\psi_0\|_s\le R$, the solution exists on $[0,\tau_R]$
and~satisfies
$$
\sup_{0\le t\le\tau_R}\|\psi(t)\|_s\le2R.
$$
\end{proposition}

Since equation~\eqref{E:0.1} is not necessarily globally well posed,
for $\psi_0\in H^s$ and $T>0$, we denote by~$\Theta(\psi_0,T)$
the set of controls $u\in L^2(0,T;\R^q)$ for which the solution
with initial state $\psi_0$ exists on $[0,T]$.
For such a control, we denote by $\cR_t(\psi_0,u)$ the value
of the solution at time $t\in[0,T]$.

\begin{definition}\label{D:1.2}
Let $\psi_0,\psi_1\in H^s\cap\bS$ and $T>0$. We say that
equation~\eqref{E:0.1} is approximately controllable in $H^s$
from $\psi_0$ to $\psi_1$ at time $T$ if, for every $\e>0$,
there is $u\in\Theta(\psi_0,T)$ such that
$$
\|\cR_T(\psi_0,u)-\psi_1\|_s<\e.
$$
It is exactly controllable in $H^s$ from $\psi_0$ to $\psi_1$ at time $T$
if $\cR_T(\psi_0,u)=\psi_1$ for some $u\in\Theta(\psi_0,T)$.
We~say~that approximate or exact controllability is global on
$H^s\cap\bS$ at time $T$ if the corresponding property holds
for every pair of states in this set.
\end{definition}

Let us recall how the saturation condition of~\cite{DN} arises.
For every $\psi_0\in H^s$ and every $\varphi\in C^\infty(\T;\R)$,
\cite[Proposition~1.2]{DN} establishes the limit
$$
e^{-i\varphi/\sqrt\e}
\cR_\e\bigl(e^{i\varphi/\sqrt\e}\psi_0,0\bigr)
\longrightarrow e^{-i(\partial_x\varphi)^2}\psi_0
\quad\text{in }H^s\quad\text{as }\e\to0^+.
$$
The solution on the left exists for all sufficiently small $\e>0$.
Suppose that, for every $\psi_0\in H^s$ and~$c\in\R$,
we can approximate $e^{ic\varphi}\psi_0$ in arbitrarily small
time by a controlled solution starting from~$\psi_0$. We can then approximate the two phase multiplications
on the left, before and after the uncontrolled evolution.
The limit gives a new phase $-(\partial_x\varphi)^2$.

To describe the phases obtained by repeating this construction,
let $E\subset C^\infty(\T;\R)$ be a finite-dimensional real
vector space. Following~\cite[Definition~2.1]{DN}, let $\cF(E)$
denote the largest vector subspace contained in
$$
\left\{\theta_0-\sum_{j=1}^n(\partial_x\theta_j)^2:
n\ge1,\quad \theta_0,\dots,\theta_n\in E\right\}.
$$
Starting from the control potentials, we set
$$
\HH_0=\operatorname{span}\{Q_1,\dots,Q_q\},\qquad
\HH_{j+1}=\cF(\HH_j),\qquad j\ge0.
$$
These spaces are finite-dimensional and $\HH_j\subset\HH_{j+1}$.
We say that $\HH_0$ is saturating if the following~holds.
\begin{itemize}
\item[\hypertarget{H1}{\textup{(H$_1$)}}]
The union $\bigcup_{j\ge0}\HH_j$ is dense in $C^{s+2}(\T;\R)$.
\end{itemize}
In particular, by~\cite[Proposition~2.6]{DN}, the subspace $\operatorname{span}\{\mathbf1,\cos x,\sin x\}$
is saturating for every integer~$s\ge1$. Thus~\eqref{E:0.3} implies
\textup{(\hyperlink{H1}{H$_1$})}.

Let us recall the phase-control theorem in the form used below.
This formulation follows from~\cite[Theorem~2.2]{DN} by approximating
the phase in $H^s$ by smooth real-valued functions.

\begin{theorem}\label{T:1.3}
Assume \textup{(\hyperlink{H1}{H$_1$})}.
For every $\psi_0\in H^s$, every $\theta\in H^s(\T;\R)$, and
every $\e,\delta>0$, there are $t\in(0,\delta)$ and
$u\in\Theta(\psi_0,t)$ such that
\begin{equation}\label{E:1.2}
\|\cR_t(\psi_0,u)-e^{i\theta}\psi_0\|_s<\e.
\end{equation}
\end{theorem}
The initial state in this theorem may have zeros. In that case, the
target $e^{i\theta}\psi_0$ vanishes at exactly the same points as
$\psi_0$. On the other hand, when the initial state does not vanish, any
target in~$H^s\cap\bS$ can be approximated. This is proved
in~\cite[Theorem~3.1]{NLSI} for nowhere vanishing targets. Since~nowhere~vanishing states are dense in $H^s\cap\bS$ on the circle,
the result extends to arbitrary targets;
see~\cite[Remark~3.3]{NLSI}.

\begin{theorem}\label{T:1.4}
Assume \textup{(\hyperlink{H1}{H$_1$})}.
For any $T>0$, any nowhere vanishing initial state
$\psi_0\in H^s\cap\bS$, and any target state
$\psi_1\in H^s\cap\bS$, equation~\eqref{E:0.1}
is approximately controllable in $H^s$ from $\psi_0$
to $\psi_1$ at time $T$.
\end{theorem}

\section{Global approximate controllability}\label{S:2}

We now extend Theorem~\ref{T:1.4} to initial states
with zeros, keeping the regularity assumptions of the previous
section on $s$ and the potentials.
The following result generalizes Theorem~\ref{T:A} by
replacing~\eqref{E:0.3} with the saturation condition~\textup{(\hyperlink{H1}{H$_1$})}.

\begin{theorem}\label{T:2.1}
Assume \textup{(\hyperlink{H1}{H$_1$})}.
For any $T>0$ and any initial and target states
$\psi_0,\psi_1\in H^s\cap\bS$, equation~\eqref{E:0.1}
is approximately controllable in $H^s$ from $\psi_0$
to $\psi_1$ at time $T$.
\end{theorem}

In Section~\ref{SS:2.1}, we first approximate the free
Schr\"odinger evolution by controlled solutions of~\eqref{E:0.1}.
In Section~\ref{SS:2.2}, we combine this approximation
with phase control to steer any nonzero initial state to a
nowhere vanishing state in arbitrarily small time.  We complete the proof in Section~\ref{SS:2.3} by
applying Theorem~\ref{T:1.4} on the remaining time interval. 
For the rest of this section, we assume
\textup{(\hyperlink{H1}{H$_1$})}.

\subsection{Approximation of the free evolution}\label{SS:2.1}

We first show that the free Schr\"odinger group $U(t)$ can be
approximated by controlled solutions of~\eqref{E:0.1}.

\begin{proposition}\label{P:2.2}
For any $\psi_0\in H^s$, $T>0$, and $\e>0$,
there is a control $u\in\Theta(\psi_0,T)$ such that
\begin{equation}\label{E:2.1}
\|\cR_T(\psi_0,u)-U(T)\psi_0\|_s<\e.
\end{equation}
\end{proposition}

Before proving the proposition, let us see how a phase change
after a short time can compensate for the potential and the
nonlinearity. We denote by $S(t)$ the local flow of
equation~\eqref{E:0.1} with $u=0$ and write the uncontrolled
equation in the form
$$
\partial_t\psi=i\partial_x^2\psi+F(\psi),
$$
where $F:H^s\to H^s$ is defined by
\begin{equation}\label{E:2.2}
F(v)=-i\bigl(V+\kappa|v|^{2p}\bigr)v,\qquad v\in H^s.
\end{equation}
Omitting the Laplacian, we obtain the auxiliary ODE
in $H^s$
$$
\partial_t w=F(w),\qquad w(0)=v.
$$
Its unique solution belongs to $C^1(\R;H^s)$
and is given by the explicit formula
$w(t)=G(t)v$, where
\begin{equation}\label{E:2.3}
(G(t)v)(x)=e^{-it(V(x)+\kappa|v(x)|^{2p})}v(x),
\qquad t\in\R.
\end{equation}
In particular,
$|(G(t)v)(x)|=|v(x)|$ for every $t\in\R$ and $x\in\T$.

To approximate $U(h)v$, we write $i\partial_x^2$ as the sum of
$i\partial_x^2+F$ and $-F$. Following the idea of the Lie--Trotter product formula for flows
(e.g., see~\cite[p.~125]{Mar74}), we let the uncontrolled equation
evolve for time $h$ and then apply $G(-h)$ to the resulting state. Thus, for $h>0$ sufficiently
small, we set
\begin{equation}\label{E:2.4}
P_hv=G(-h)S(h)v
=e^{ih(V+\kappa|S(h)v|^{2p})}S(h)v.
\end{equation}
Figure~\ref{F:free-approximation} illustrates this step.

\par\bigskip\medskip
\noindent\begin{minipage}{\linewidth}
\captionsetup{type=figure}
\centering
\captionsetup{skip=7pt}
\scalebox{1.06}{%
\begin{tikzpicture}[x=1cm,y=1.30cm,
  every node/.style={font=\small,inner sep=1.5pt},
  drift/.style={draw=black!83,line width=.92pt},
  kick/.style={draw=darkblue,line width=1.22pt},
  orbit/.style={draw=teal!64!black,line width=.95pt,
                dash pattern=on 3.4pt off 2.1pt},
  separation/.style={draw=black!58,line width=.58pt,densely dotted},
  direction/.style={postaction={decorate},
    decoration={markings,mark=at position #1 with
      {\arrow{Stealth[length=4.8pt,width=3.9pt]}}}},
  state/.style={circle,fill=darkblue,draw=white,line width=.45pt,inner sep=2pt},
  freepoint/.style={circle,fill=white,draw=teal!64!black,
                    line width=.85pt,inner sep=1.7pt}]
\coordinate (start)  at (0,0);
\coordinate (driftend) at (4.65,2.02);
\coordinate (splitend) at (7.10,.95);
\coordinate (freeend) at (7.24,.20);

\draw[separation] (splitend) -- (freeend);
\draw[orbit,direction=.57] (start)
  .. controls (2.01,-.31) and (5.35,-.13) .. (freeend);
\draw[drift,direction=.53] (start)
  .. controls (1.59,.23) and (2.11,1.94) .. (driftend);
\draw[kick,direction=.55] (driftend)
  .. controls (6.03,2.06) and (6.54,1.39) .. (splitend);

\fill[black!85] (start) circle (2pt);
\node[below left] at (start) {$v$};
\fill[black!82] (driftend) circle (1.45pt);
\node[above=2pt] at (driftend) {$S(h)v$};
\node[state] at (splitend) {};
\node[above right=2pt] at (splitend) {$P_hv$};
\node[freepoint] at (freeend) {};
\node[below right=2pt] at (freeend) {$U(h)v$};

\node[above left=2pt] at (2.17,1.09) {$S(r)v$};
\node[above right=3pt,text=darkblue] at (6.16,1.73) {$G(-r)S(h)v$};
\node[below=3pt] at (3.56,-.15) {$U(r)v$};
\draw[black!45,line width=.4pt] (7.25,.57) -- (7.71,.57);
\node[anchor=west] at (7.76,.57)
  {$\|P_hv-U(h)v\|_s\le h\varrho(h)$};

\begin{scope}[yshift=-.85cm]
  \draw[drift] (.15,0) -- (.78,0);
  \node[right] at (.87,0) {zero control};
  \draw[kick,direction=.55] (3.42,0) -- (4.05,0);
  \node[right] at (4.14,0) {exact phase change};
  \draw[orbit] (7.75,0) -- (8.38,0);
  \node[right] at (8.47,0) {free solution};
\end{scope}
\end{tikzpicture}%
}
\caption{Starting from $v$ on the free trajectory, we follow the uncontrolled flow
to $S(h)v$ and apply the phase change $G(-h)$ to reach $P_hv$.
The dotted segment shows the error in~\eqref{E:2.6}.}
\label{F:free-approximation}
\end{minipage}
\par\bigskip\medskip

To iterate $P_h$, we need to estimate its dependence on the initial state and the difference between $P_hv$ and~$U(h)v$. For $R>0$, let~us~set
$$
B_R=\{v\in H^s:\|v\|_s\le R\}.
$$
\begin{lemma}\label{L:2.3}
For any $\psi_0\in H^s$ and $T>0$, let
$$
R=\|\psi_0\|_s+2,\qquad
K=\{U(t)\psi_0:0\le t\le T\}\subset B_R.
$$
There are constants $h_0,C>0$ and a function
$\varrho:(0,h_0]\to[0,\infty)$, with $\varrho(h)\to0$ as
$h\to0^+$, such that the following properties hold for
every $h\in(0,h_0]$.
\begin{enumerate}
\item[\textup{(i)}]
For $0\le r\le h$, the flow $S(r)$ is defined on $B_R$ and
$$
S(r)B_R\subset B_{2R},\qquad G(-r)B_{2R}\subset B_{3R}.
$$
 Moreover, 
\begin{equation}\label{E:2.5}
\|P_hv-P_hw\|_s\le e^{Ch}\|v-w\|_s, \quad v,w\in B_R.
\end{equation}
\item[\textup{(ii)}]
We have
\begin{equation}\label{E:2.6}
\sup_{v\in K}\|P_hv-U(h)v\|_s\le h\varrho(h).
\end{equation}
\end{enumerate}
\end{lemma}

\begin{proof}
Since $H^s$ is an algebra, $M_R=\sup_{v\in B_{3R}}\|F(v)\|_s$ is finite
and $F$ is Lipschitz on $B_{3R}$.
Let $L_R>0$ be a corresponding Lipschitz constant.

By Proposition~\ref{P:1.1}, there is $h_0>0$
such that $S(r)$ is defined on $B_R$ and
$S(r)B_R\subset B_{2R}$ for $0\le r\le h_0$.
For the auxiliary flow, we have
\begin{equation}\label{E:2.7}
G(-r)z=z-\int_0^r F(G(-a)z)\,\dd a,
\qquad z\in H^s,\quad r\ge0.
\end{equation}
If $z\in B_{2R}$, then, as long as $G(-r)z$
remains in $B_{3R}$, \eqref{E:2.7} implies
$$
\|G(-r)z\|_s\le2R+rM_R.
$$
Decreasing $h_0$ if necessary so that $h_0M_R<R$, we see
by continuity that the trajectory $r\mapsto G(-r)z$ cannot
leave $B_{3R}$ on $[0,h_0]$. 
Thus both inclusions in (i) hold.

Applying Gronwall's inequality to the Duhamel formula
for $S(r)$ and to~\eqref{E:2.7}, we obtain
\begin{align*}
\|S(r)v-S(r)w\|_s
&\le e^{L_Rr}\|v-w\|_s,
&&v,w\in B_R,\\
\|G(-r)z-G(-r)\widetilde z\|_s
&\le e^{L_Rr}\|z-\widetilde z\|_s,
&&z,\widetilde z\in B_{2R}
\end{align*}
for $0\le r\le h_0$.
Since $P_h=G(-h)S(h)$, these estimates
imply~\eqref{E:2.5} with $C=2L_R$.
This proves (i).

To prove (ii), let $0<h\le h_0$ and $v\in K$.
Since $P_hv=G(-h)S(h)v$, taking $S(h)v$ as the initial state
in~\eqref{E:2.7} gives
$$
P_hv-S(h)v=-\int_0^h F(G(-r)S(h)v)\,\dd r.
$$
On the other hand, by Duhamel's formula for $S(h)v$,
$$
S(h)v-U(h)v=\int_0^h U(h-r)F(S(r)v)\,\dd r.
$$
Adding these two identities, we obtain
\begin{equation}\label{E:2.8}
P_hv-U(h)v
=\int_0^h\bigl[U(h-r)F(S(r)v)-F(G(-r)S(h)v)\bigr]\,\dd r.
\end{equation}

Let $D(h,r,v)$ denote the integrand in~\eqref{E:2.8}
and set
$$
\Delta=\{(h,r)\in\R^2:0\le r\le h\le h_0\}.
$$
By Proposition~\ref{P:1.1} and~\eqref{E:2.3}, the flows $S$ and $G$ depend continuously on time and on the
initial state. The~same holds for $U$ by strong continuity and
the isometry of $U(t)$ in $H^s$. Since $F$ is continuous,
the map
$$
D:\Delta\times K\longrightarrow H^s
$$
is continuous. Moreover, $\Delta$ is compact in $\R^2$, and
$K$ is compact in $H^s$ as the image of $[0,T]$ under the
continuous map $t\mapsto U(t)\psi_0$. Thus $D$ is uniformly
continuous on $\Delta\times K$.

At $h=r=0$, the two terms defining $D(h,r,v)$ are both equal
to $F(v)$. Thus $D(0,0,v)=0$ for every $v\in K$, and the
uniform continuity of $D$ on $\Delta\times K$ gives
$$
\varrho(h):=\sup_{\substack{v\in K\\0\le r\le h}}
\|D(h,r,v)\|_s\longrightarrow0
\qquad\text{as }h\to0^+.
$$
Now~\eqref{E:2.6} follows from~\eqref{E:2.8}.
This proves~(ii).
\end{proof}

\begin{proof}[Proof of Proposition~\ref{P:2.2}]
Let us fix $\psi_0\in H^s$ and $T>0$, and use the notation of
Lemma~\ref{L:2.3}. At each step, we first use zero
control for time $h$ and then apply a short control to approximate
the phase change in~\eqref{E:2.4}. To allow time for the phase
controls, we take $h=T/(n+1)$, where $n$ is the number of steps.
The $n$ uncontrolled intervals then have total length $T-h$.
We choose the phase controls so that their total duration is
less than $h$ and use zero control for the remaining time.

\emph{Step~1. Construction of the control.}
Let us take $n$ large enough, so that $h\le h_0$. We construct the control
for equation~\eqref{E:0.1} starting from $g_0=\psi_0$.
Suppose that, after the first $j$ steps, the controlled solution
is at $g_j\in B_R$, with $0\le j<n$. We use zero control for
the next interval of length $h$ and denote by $b_j=S(h)g_j$
the state at the end of this interval. By~\eqref{E:2.4},
$$
P_hg_j=e^{i\theta_j}b_j,\qquad
\theta_j=h\bigl(V+\kappa|b_j|^{2p}\bigr)\in H^s(\T;\R).
$$
Starting from $b_j$, we now use Theorem~\ref{T:1.3} to
approximate $P_hg_j$. There are $\sigma_j\in(0,h/n)$ and
$u_j\in\Theta(b_j,\sigma_j)$ such that
$$
g_{j+1}=\cR_{\sigma_j}(b_j,u_j),\qquad
\|g_{j+1}-P_hg_j\|_s<h^2.
$$
Thus we reach $g_{j+1}$ from $g_j$ in time $h+\sigma_j$.
To see that this construction can be continued up to $g_n$,
let us compare $g_j$ with the free solution $z_j=U(jh)\psi_0$.
Here $jh$ is the total length of the first $j$ intervals on
which we use zero control. Let~$E_j=\|g_j-z_j\|_s$.
Since $z_j\in K$ and $z_{j+1}=U(h)z_j$, estimates
\eqref{E:2.5} and~\eqref{E:2.6} give, as long
as $g_j\in B_R$,
\begin{align*}
E_{j+1}
&\le \|g_{j+1}-P_hg_j\|_s+\|P_hg_j-P_hz_j\|_s
     +\|P_hz_j-U(h)z_j\|_s\\
&\le h^2+e^{Ch}E_j+h\varrho(h).
\end{align*}
Since $E_0=0$ and $nh<T$, iteration gives
\begin{equation}\label{E:2.9}
E_k\le h\bigl(h+\varrho(h)\bigr)\sum_{j=0}^{k-1}e^{Chj}
\le Te^{CT}\bigl(h+\varrho(h)\bigr)
\end{equation}
for every $k\le n$ for which $g_k$ has been constructed.
For all sufficiently large $n$, the last term is less than~$1$.
Since $\|z_k\|_s=\|\psi_0\|_s$, we then have
$$
\|g_k\|_s\le\|\psi_0\|_s+E_k<\|\psi_0\|_s+1<R.
$$
It follows by induction that all the states $g_0,\ldots,g_n$
are defined and belong to $B_R$.

\emph{Step~2. Completion at time $T$.}
After $n$ steps, the controlled solution is at $g_n$ at time
$nh+\sigma$, where~$\sigma=\sum_{j=0}^{n-1}\sigma_j$.
Since each $\sigma_j<h/n$, we have $\sigma<h$. We use zero
control for the remaining time
$$
r=T-(nh+\sigma)=h-\sigma>0.
$$
Figure~\ref{F:control-times} illustrates the construction for $n=3$.

\par\bigskip\medskip
\noindent\begin{minipage}{\linewidth}
\centering
\captionsetup{type=figure}
\usetikzlibrary{arrows.meta,decorations.markings}
\captionsetup{skip=7pt}
\scalebox{1.06}{%
\begin{tikzpicture}[x=1cm,y=1.30cm,
  every node/.style={font=\small,inner sep=1.5pt},
  drift/.style={draw=black!83,line width=.92pt},
  kick/.style={draw=darkblue,line width=1.22pt},
  orbit/.style={draw=teal!64!black,line width=.95pt,
                dash pattern=on 3.4pt off 2.1pt},
  separation/.style={draw=black!58,line width=.58pt,densely dotted},
  direction/.style={postaction={decorate},
    decoration={markings,mark=at position #1 with
      {\arrow{Stealth[length=4.8pt,width=3.9pt]}}}},
  state/.style={circle,fill=darkblue,draw=white,line width=.45pt,inner sep=2pt},
  freepoint/.style={circle,fill=white,draw=teal!64!black,
                    line width=.85pt,inner sep=1.7pt}]
\path (0,3.65) -- (0,3.64);
\coordinate (start) at (0,0);
\coordinate (b0) at (2.06,1.78);
\coordinate (b1) at (5.14,2.70);
\coordinate (b2) at (8.32,2.83);
\coordinate (g1) at (3.05,.89);
\coordinate (g2) at (6.23,1.70);
\coordinate (g3) at (9.52,1.28);
\coordinate (z1) at (3.19,.26);
\coordinate (z2) at (6.37,1.07);
\coordinate (z3) at (9.64,.56);
\coordinate (nearend) at (11.32,.37);
\coordinate (target) at (13.12,.70);
\coordinate (finish) at (12.95,1.31);

\foreach \j in {1,2,3}{\draw[separation] (g\j) -- (z\j);}
\draw[separation] (finish) -- (target);

\draw[orbit,direction=.57] (start) .. controls (.96,-.17) and (2.29,-.01) .. (z1);
\draw[orbit,direction=.55] (z1) .. controls (4.17,.56) and (5.46,1.04) .. (z2);
\draw[orbit,direction=.56] (z2) .. controls (7.55,1.11) and (8.77,.73) .. (z3);
\draw[orbit] (z3) .. controls (10.30,.43) and (10.83,.36) .. (nearend);
\draw[orbit,direction=.60] (nearend) .. controls (12.01,.38) and (12.55,.54) .. (target);

\draw[drift,direction=.53] (start) .. controls (.63,.19) and (.87,1.72) .. (b0);
\draw[kick,direction=.52] (b0) .. controls (2.82,1.83) and (2.82,1.19) .. (g1);
\draw[drift,direction=.53] (g1) .. controls (3.89,.98) and (3.99,2.67) .. (b1);
\draw[kick,direction=.52] (b1) .. controls (5.90,2.74) and (5.99,2.09) .. (g2);
\draw[drift,direction=.52] (g2) .. controls (7.09,1.86) and (7.30,2.98) .. (b2);
\draw[kick,direction=.56] (b2) .. controls (9.20,2.70) and (9.10,1.78) .. (g3);
\draw[drift,direction=.55] (g3) .. controls (10.43,1.40) and (11.68,1.66) .. (finish);

\fill[black!85] (start) circle (2pt);
\node[below left] at (start) {$\psi_0$};
\foreach \j/\prev in {1/0,2/1,3/2}{
  \fill[black!82] (b\prev) circle (1.45pt);
  \node[above=2pt] at (b\prev) {$b_{\prev}$};
  \node[state] at (g\j) {};
  \node[below right=1pt] at (g\j) {$g_{\j}$};
  \node[freepoint] at (z\j) {};
  \node[below=2pt] at (z\j) {$z_{\j}$};
}
\node[above left=1pt] at (.94,1.12) {$S(h)$};
\node[right=2pt,text=darkblue] at (2.80,1.45) {$u_0$};
\node[left=3pt] at (6.30,1.40) {$E_2$};
\node[above=2pt] at (11.08,1.50) {$S(r)$};
\node[state] at (finish) {};
\node[above=2pt] at (finish) {$S(r)g_3$};
\node[freepoint] at (nearend) {};
\node[below=3pt] at (nearend) {$U(T-\sigma)\psi_0$};
\node[freepoint] at (target) {};
\node[below right=2pt] at (target) {$U(T)\psi_0$};

\begin{scope}[yshift=-.89cm]
  \draw[drift] (1.05,0) -- (1.68,0);
  \node[right] at (1.77,0) {zero control};
  \draw[kick,direction=.55] (4.47,0) -- (5.10,0);
  \node[right] at (5.19,0) {phase control};
  \draw[orbit] (8.17,0) -- (8.80,0);
  \node[right] at (8.89,0) {free solution};
\end{scope}
\end{tikzpicture}%
}
\caption{The construction for $n=3$. From $g_j$, we follow the
uncontrolled flow to $b_j=S(h)g_j$ and apply $u_j$ for time
$\sigma_j$ to reach $g_{j+1}$. The dotted segments join $g_j$
to $z_j=U(jh)\psi_0$ and $S(r)g_3$ to   $U(T)\psi_0$.}
\label{F:control-times}
\end{minipage}
\par\bigskip\medskip

Since $g_n\in B_R$ and $r<h\le h_0$, part~\textup{(i)} of
Lemma~\ref{L:2.3} gives $S(t)g_n\in B_{2R}$ for
$0\le t\le r$. Thus the concatenated control $u^{(n)}$ belongs
to $\Theta(\psi_0,T)$ and
$\cR_T(\psi_0,u^{(n)})=S(r)g_n$.

The uncontrolled intervals have total length $nh+r=T-\sigma$.
Since $z_n=U(nh)\psi_0$, the free solution at this time is
$U(T-\sigma)\psi_0=U(r)z_n$. Subtracting $U(r)z_n$ from
Duhamel's formula for $S(r)g_n$, we obtain
$$
S(r)g_n-U(T-\sigma)\psi_0
=U(r)(g_n-z_n)+\int_0^r U(r-a)F(S(a)g_n)\,\dd a.
$$
Since $U$ is an isometry in $H^s$ and $S(a)g_n\in B_{2R}$
for $0\le a\le r$, it follows that
\begin{equation}\label{E:2.10}
\|S(r)g_n-U(T-\sigma)\psi_0\|_s
\le E_n+rM_R,\qquad M_R=\sup_{v\in B_{3R}}\|F(v)\|_s.
\end{equation}

By~\eqref{E:2.9}, $E_n\to0$ as $n\to\infty$.
Moreover, both $r$ and $\sigma$ tend to zero, since
$0<r,\sigma<h=T/(n+1)$. Thus the right-hand side
of~\eqref{E:2.10} tends to zero, while
$$
U(T-\sigma)\psi_0\to U(T)\psi_0 \quad \text{in} \,\,\, H^s
$$ by strong
continuity of $U$. We conclude that the final state
$\cR_T(\psi_0,u^{(n)})=S(r)g_n$ converges to $U(T)\psi_0$
in $H^s$, which proves~\eqref{E:2.1}.
\end{proof}

\subsection{Removal of zeros}
\label{SS:2.2}

We now show how to steer any 
nonzero initial state $\psi_0$ 
to a nowhere vanishing state in arbitrarily small time using the
controls in~\eqref{E:0.1}. For this, we first choose a periodic
phase $\theta$ and a short time $\tau>0$ such that
$U(\tau)(e^{i\theta}\psi_0)$ is nowhere vanishing.
\begin{lemma}\label{L:2.4}
For any $\psi_0\in H^s\setminus\{0\}$ and $a>0$,
there are $\tau\in(0,a)$ and~$\theta\in C^\infty(\T;\R)$ such that~$U(\tau)(e^{i\theta}\psi_0)$ is nowhere vanishing.
\end{lemma}

\begin{proof}
Since $\psi_0$ is continuous and not identically zero,
we can choose an open interval $I\subset\{\psi_0\ne0\}$
with $0<|I|<2\pi$.
Let $N\ge3$ be an odd integer large enough that
$$
\tau=\frac{2\pi}{N}<\min\{a,|I|/3\}.
$$
Then, for every $x\in\T$, at least three points
of the grid $\{x+k\tau:0\le k<N\}$ belong to $I$. For this choice of $\tau$, the Fourier multiplier $e^{-i\tau n^2}$
is $N$-periodic in $n$.
Its discrete Fourier expansion gives, for~$v\in H^s$,
\begin{equation}\label{E:2.11}
U(\tau)v(x)=\sum_{k=0}^{N-1}c_kv(x+k\tau),
\qquad
c_k=\frac1N\sum_{j=0}^{N-1}e^{-i\tau(j^2+kj)}.
\end{equation}
Indeed, for every $n\in\Z$,
$$
\sum_{k=0}^{N-1}c_ke^{ink\tau}
=\frac1N\sum_{j=0}^{N-1}e^{-i\tau j^2}
\sum_{k=0}^{N-1}e^{i\tau(n-j)k}
=e^{-i\tau n^2}.
$$
Thus the two sides of~\eqref{E:2.11} have the same
Fourier coefficients, which proves the identity.

Let us check that all the coefficients $c_k$
are nonzero.
Since $N$ is odd, the change of variables
$\alpha=j-\ell$, $\beta=j+\ell$ is invertible
modulo $N$. Hence
$$
N^2|c_k|^2
=\sum_{j,\ell=0}^{N-1}
e^{-i\tau(j-\ell)(j+\ell+k)}
=\sum_{\alpha,\beta=0}^{N-1}
e^{-i\tau\alpha(\beta+k)}=N.
$$
It follows that $|c_k|=N^{-1/2}$ for every $k$.

Let us now set $f_k(x)=c_k\psi_0(x+k\tau)$.
Since $s\ge1$, each $f_k$ belongs to
$H^1(\T;\C)$.
By the choice of $\tau$ and the fact that $c_k\ne0$,
at least three of the values $f_k(x)$ are nonzero
at every $x\in\T$.
Applying Lemma~\ref{L:A.1}, we find
$\theta\in C^\infty(\T;\R)$ such that
$$
U(\tau)(e^{i\theta}\psi_0)(x)
=\sum_{k=0}^{N-1}f_k(x)e^{i\theta(x+k\tau)}\ne0,
\qquad x\in\T.
$$
This proves the lemma.
\end{proof}

We now return to equation~\eqref{E:0.1} and use phase control
and Proposition~\ref{P:2.2} to approximate the nowhere
vanishing state given by Lemma~\ref{L:2.4}.

\begin{lemma}\label{L:2.5}
For any $\psi_0\in H^s\setminus\{0\}$ and $\delta>0$,
there are $t\in(0,\delta)$ and
$u\in\Theta(\psi_0,t)$ such that
$\cR_t(\psi_0,u)$ is nowhere vanishing.
\end{lemma}
\begin{proof}
By Lemma~\ref{L:2.4}, we can choose
$\tau\in(0,\delta/2)$ and $\theta\in C^\infty(\T;\R)$
such that $h=U(\tau)(e^{i\theta}\psi_0)$ is nowhere vanishing.
Let us set
\begin{equation}\label{E:2.12}
m=\min_{x\in\T}|h(x)|>0.
\end{equation}
 Let $C_s>0$ be such that
$\|v\|_{L^\infty}\le C_s\|v\|_s$ for all $v\in H^s$.
By Theorem~\ref{T:1.3}  there are $\sigma\in(0,\delta/2)$
and~$u_0\in\Theta(\psi_0,\sigma)$ such that
\begin{equation}\label{E:2.13}
g=\cR_\sigma(\psi_0,u_0),\qquad
\|g-e^{i\theta}\psi_0\|_s<\frac{m}{4C_s}.
\end{equation}
Starting from $g$, we apply Proposition~\ref{P:2.2}
at time $\tau$ to choose $u_1\in\Theta(g,\tau)$ such that
\begin{equation}\label{E:2.14}
\|\cR_\tau(g,u_1)-U(\tau)g\|_s<\frac{m}{4C_s}.
\end{equation}
Let $u$ be the concatenation of $u_0$ and $u_1$, and set
$t=\sigma+\tau<\delta$. Since $\cR_\sigma(\psi_0,u_0)=g$,
we have $u\in\Theta(\psi_0,t)$ and
$\cR_t(\psi_0,u)=\cR_\tau(g,u_1)$.
Using~\eqref{E:2.13}, \eqref{E:2.14},
and the isometry of $U(\tau)$ in $H^s$, we obtain
$$
\|\cR_t(\psi_0,u)-h\|_{L^\infty}
\le C_s\bigl(\|\cR_\tau(g,u_1)-U(\tau)g\|_s
             +\|g-e^{i\theta}\psi_0\|_s\bigr)
<\frac m2.
$$
Hence, for every $x\in\T$,
$$
|\cR_t(\psi_0,u)(x)|
\ge |h(x)|-\|\cR_t(\psi_0,u)-h\|_{L^\infty}
>\frac m2.
$$
\end{proof}

\subsection{Proof of Theorem~\ref{T:2.1}}
\label{SS:2.3}

By Lemma~\ref{L:2.5}, we can choose $t_0\in(0,T)$
and $u_0\in\Theta(\psi_0,t_0)$ such that the state
$h_0=\cR_{t_0}(\psi_0,u_0)$ is nowhere vanishing.
Applying Theorem~\ref{T:1.4} with initial
state $h_0$ and time $T-t_0$, we find a control~$u_1\in\Theta(h_0,T-t_0)$ such~that
$$
\|\cR_{T-t_0}(h_0,u_1)-\psi_1\|_s<\e.
$$
Concatenating $u_0$ and $u_1$, we obtain
$u\in\Theta(\psi_0,T)$ satisfying
$$
\|\cR_T(\psi_0,u)-\psi_1\|_s<\e.
$$
This completes the proof of the theorem.

\section{Well-posedness for the additional potentials}\label{S:3}

Here we prove local well-posedness in $H^2$ for equation~\eqref{E:0.5}
with $Q=(\mathbf1,\cos x,\sin x,Q_c,Q_s)$, where $Q_c$ and~$Q_s$
are given by~\eqref{E:0.4}. Since multiplication
by these two potentials does not preserve $H^2$, we estimate
the corresponding Duhamel terms in $H^2$, following
\cite[Sections~3.1--3.2]{BL}.

For $T>0$, let
$$
X_T=C([0,T];H^2),\qquad
\|\psi\|_{X_T}=\sup_{0\le t\le T}\|\psi(t)\|_2.
$$
By a solution of~\eqref{E:0.5} on $[0,T]$, we mean a function
$\psi\in X_T$ satisfying the Duhamel formula in $H^1$.
\subsection{Duhamel smoothing estimate}

The potentials $Q_c$ and $Q_s$ are smooth away
from $x=0$. For $P\in\{Q_c,Q_s\}$ and $\psi\in H^2$,
we write
$$
P\psi=P(\psi-\psi(0))+\psi(0)P.
$$
Since $\psi-\psi(0)$ vanishes at $x=0$, the following lemma
 shows that the first term belongs
to $H^2$.
The second term need not belong to $H^2$, but its
Duhamel integral can be estimated in this space;
see Lemma~\ref{L:3.3}.

\begin{lemma}\label{L:3.1}
For $P\in\{Q_c,Q_s\}$ and $\psi\in H^2$, let
$$
R_P\psi(x)=P(x)\bigl(\psi(x)-\psi(0)\bigr).
$$
Then $R_P$ is a bounded linear operator on $H^2$.
\end{lemma}

\begin{proof}
On $[-\pi,\pi]$, the series defining $Q_c$ and $Q_s$
have the representations
$$
Q_c(x)=\frac{\pi^2}{6}-\frac{\pi|x|}{2}+\frac{x^2}{4},
\qquad
Q_s(x)=-\int_0^x\log\bigl(2|\sin(y/2)|\bigr)\,\dd y.
$$
Thus, for $0<|x|<\pi$,
$$
Q_c'(x)=\frac{x}{2}-\frac{\pi}{2}\sgn(x),
\qquad
Q_s'(x)=-\log\bigl(2|\sin(x/2)|\bigr),
$$
and
$$
Q_c''(x)=\frac12,\qquad
Q_s''(x)=-\frac12\cot(x/2).
$$
For either choice of $P$, these formulas give
$P\in L^\infty$, $P'\in L^2$, and
$|P''(x)|\le C/|x|$ for $0<|x|<\pi$.

Let $v=\psi-\psi(0)$. By Sobolev embedding,
$$
\|v\|_2+\|v'\|_{L^\infty}\le C\|\psi\|_2.
$$
Since $v(0)=0$, we have
$|v(x)|\le |x|\,\|v'\|_{L^\infty}$, and hence
$$
\|P''v\|_{L^2}\le C\|v'\|_{L^\infty}.
$$
Moreover, the formulas for $P'$ show that
$P'(x)v(x)\to0$ as $x\to0^\pm$.
Since $P$ and $v'$ are continuous, the first derivative
$$
(Pv)'=P'v+Pv'
$$
has matching one-sided limits at zero.
The boundary terms at zero therefore cancel, and
integration by parts~gives
$$
(Pv)''=P''v+2P'v'+Pv''
$$
in the sense of distributions on $\T$. Here $P''$ on the right-hand side denotes the
classical derivative away from zero.
All three terms belong to $L^2$, and the preceding
estimates yield
$$
\|Pv\|_2
\le C\bigl(\|v\|_2+\|v'\|_{L^\infty}\bigr)
\le C\|\psi\|_2.
$$
This proves the lemma.
\end{proof}

To estimate the Duhamel integral corresponding to the term $\psi(0)P$,
we use the following consequence of Bessel's~inequality.
\begin{lemma}\label{L:3.2}
For any $T>0$, there is $C_T>0$ such that, for any
interval $I\subset[0,T]$ and $f\in L^2(I;\C)$, we~have
$$
\sum_{n\in\Z}
\Big|\int_I e^{in^2r}f(r)\,\dd r\Big|^2
\le C_T\|f\|_{L^2(I)}^2.
$$
\end{lemma}

\begin{proof}
Let us take an integer $N\ge1$ such that $T\le2\pi N$
and extend $f$ by zero to $(0,2\pi N)$.
Applying Bessel's inequality to the orthonormal family
$\{(2\pi N)^{-1/2}e^{-ikr}:k\in\Z\}$ in
$L^2(0,2\pi N)$, we obtain the
required estimate with $C_T=4\pi N$.
\end{proof}
For $0\le t_0<t_1$, $u\in L^2(t_0,t_1;\R^5)$,
and $h\in C([t_0,t_1];H^2)$, let us set
$$
B_{t_0}(u,h)(t)
=\int_{t_0}^t U(t-r)
\lan u(r),Q\ran h(r)\,\dd r,
\qquad t\in[t_0,t_1].
$$
This integral is well defined in $H^1$, since
$Q_j\in H^1$ and $H^1(\T)$ is an algebra.
We prove the following smoothing estimate by adapting
the argument of~\cite[Lemma~3]{BL}.

\begin{lemma}\label{L:3.3}
For any $T>0$, there is $C_T>0$ such that, for any
$0\le t_0<t_1\le T$, $u\in L^2(t_0,t_1;\R^5)$,
and~$h\in C([t_0,t_1];H^2)$, we have
$B_{t_0}(u,h)\in C([t_0,t_1];H^2)$ and
$$
\sup_{t\in[t_0,t_1]}\|B_{t_0}(u,h)(t)\|_2
\le C_T\|u\|_{L^2(t_0,t_1)}
\sup_{t\in[t_0,t_1]}\|h(t)\|_2.
$$
\end{lemma}

\begin{proof}
For $j=1,2,3$, multiplication by $Q_j$ is bounded
on $H^2$, so
$$
\|u_jQ_jh\|_{L^1(t_0,t_1;H^2)}
\le C_T\|u_j\|_{L^2(t_0,t_1)}
\sup_{t\in[t_0,t_1]}\|h(t)\|_2.
$$
The corresponding Duhamel integral belongs to
$C([t_0,t_1];H^2)$ and is bounded by the same quantity.
For~$j=4,5$, we write
$$
Q_jh(r)=R_{Q_j}h(r)+h(r,0)Q_j.
$$
By Lemma~\ref{L:3.1}, $R_{Q_j}$ is bounded on $H^2$,
so the preceding argument also applies to the
Duhamel integral corresponding to $u_jR_{Q_j}h$.

It remains to estimate integrals of the form
$$
G(t)=\int_{t_0}^t U(t-r)f(r)P\,\dd r,
\qquad P\in\{Q_c,Q_s\},\quad f\in L^2(t_0,t_1;\C).
$$
The series defining $Q_c$ and $Q_s$ give
\begin{equation}\label{E:3.1}
\hat Q_c(n)=\frac1{2n^2},\qquad
\hat Q_s(n)=-\frac{i\,\sgn(n)}{2n^2},\quad n\ne0,
\qquad \hat Q_c(0)=\hat Q_s(0)=0.
\end{equation}
Hence $\widehat{G(t)}(0)=0$ and
$$
|\widehat{G(t)}(n)|
=\frac1{2n^2}
\Big|\int_{t_0}^t e^{in^2r}f(r)\,\dd r\Big|,
\qquad n\ne0.
$$
By the definition of the $H^2$ norm and
Lemma~\ref{L:3.2}, we obtain
\begin{equation}\label{E:3.2}
\sup_{t\in[t_0,t_1]}\|G(t)\|_2
\le C_T\|f\|_{L^2(t_0,t_1)}.
\end{equation}
Applying this estimate on $[r,t]\subset[t_0,t_1]$
and using the unitarity of $U(t)$
on $H^2$, we get
$$
\big\|U(-t)G(t)
-U(-r)G(r)\big\|_2
 =\Big\|\int_r^t U(-\sigma)
f(\sigma)P\,\dd\sigma\Big\|_2 \le C_T\|f\|_{L^2(r,t)}.
$$
Thus $t\mapsto U(-t)G(t)$ is continuous
in $H^2$. Strong continuity of $U(t)$
then gives $G\in C([t_0,t_1];H^2)$.
Finally, for $j=4,5$, we apply
\eqref{E:3.2} with
$f(r)=u_j(r)h(r,0)$. By Sobolev embedding,
$$
\|f\|_{L^2(t_0,t_1)}
\le C\|u_j\|_{L^2(t_0,t_1)}
\sup_{t\in[t_0,t_1]}\|h(t)\|_2.
$$
Combining these estimates, we complete the proof.
\end{proof}

\subsection{The solution map}

Using Lemma~\ref{L:3.3}, we prove local well-posedness for
equation~\eqref{E:0.5} in $H^2$ and smooth dependence of the
solution on the initial state and the control.
For similar results with Neumann boundary conditions,
see~\cite[Propositions~8 and~10]{BL}.
We write $\cN(\psi)=\kappa|\psi|^{2p}\psi$.

\begin{proposition}\label{P:3.4}
Let $Q=(\mathbf1,\cos x,\sin x,Q_c,Q_s)$ and $T>0$.
\begin{itemize}
\item[(i)] For any $\psi_0\in H^2$ and $u\in L^2(0,T;\R^5)$,
there is $\tau\in(0,T]$ such that equation~\eqref{E:0.5}
with data $(\psi_0,u)$ has a unique solution in $X_\tau$.

\item[(ii)] Let $\psi^0\in X_T$ be a solution with data
$(\psi_0^0,u^0)$. There is a neighborhood $\mathcal V$ of
$(\psi_0^0,u^0)$ in $H^2\times L^2(0,T;\R^5)$
on which the solution map $\Psi:\mathcal V\to X_T$ is well
defined and $C^\infty$ in the real Fr\'echet sense.
For any~$(h_0,w)\in H^2\times L^2(0,T;\R^5)$,
the derivative $v=D\Psi(\psi_0^0,u^0)(h_0,w)$
is the unique solution in $X_T$ of the linearized equation
\begin{equation}\label{E:3.3}
i\partial_t v=-\partial_x^2v+\cN'(\psi^0)v
+\lan u^0,Q\ran v+\lan w,Q\ran\psi^0,
\qquad v(0)=h_0.
\end{equation}
Here $\cN'(\psi)$ is the real Fr\'echet derivative
of $\cN$, given by
$$
\cN'(\psi)v
=\kappa(p+1)|\psi|^{2p}v
+\kappa p|\psi|^{2p-2}\psi^2\bar v.
$$
\end{itemize}
\end{proposition}

\begin{proof}
\emph{Step~1. Local existence.}
Let us fix $T>0$, $\psi_0\in H^2$, and $u\in L^2(0,T;\R^5)$.
Take $R>2\|\psi_0\|_2$ and let $C_R$ be a Lipschitz constant
for $\cN$ on the ball of radius $R$ in $H^2$. We choose
$\tau\in(0,T]$ so small that
$$
C_R\tau+C_T\|u\|_{L^2(0,\tau)}\le\frac12.
$$
By Lemma~\ref{L:3.3}, the right-hand side of the Duhamel
formula has Lipschitz constant at most $1/2$ on the closed ball
of radius $R$ in $X_\tau$. Its value at zero is $U(\cdot)\psi_0$,
whose $X_\tau$ norm is $\|\psi_0\|_2<R/2$. Thus it maps this
ball into itself, and the contraction mapping theorem gives
a solution in $X_\tau$.

\emph{Step~2. Uniqueness.}
Let $\psi,\tilde\psi\in X_S$ be two solutions with
the same data, where $0<S\le T$, and denote~$h=\psi-\tilde\psi$.
Let us take $R>\max\{\|\psi\|_{X_S},\|\tilde\psi\|_{X_S}\}$.
On an interval $[a,b]\subset[0,S]$ with $h(a)=0$,
subtracting the Duhamel formulas and applying
Lemma~\ref{L:3.3}, we obtain
$$
\sup_{t\in[a,b]}\|h(t)\|_2
\le
\bigl(C_R(b-a)+C_T\|u\|_{L^2(a,b)}\bigr)
\sup_{t\in[a,b]}\|h(t)\|_2.
$$
By absolute continuity of the integral of $|u|^2$, the factor
on the right is at most $1/2$ whenever $b-a$ is sufficiently
small, uniformly in $a$. Since $h(0)=0$, applying the estimate
on successive intervals of this length gives $h=0$ on $[0,S]$.
Taking $S=\tau$ proves~(i).

\emph{Step~3. Smooth dependence.}
Let $\psi^0\in X_T$ be the solution with data
$(\psi_0^0,u^0)$. We define
$$
\Phi:X_T\times H^2\times L^2(0,T;\R^5)\longrightarrow X_T
$$
by
$$
\Phi(\psi,\psi_0,u)(t)
=\psi(t)-U(t)\psi_0
+i\int_0^t U(t-r)\cN(\psi(r))\,\dd r
+iB_0(u,\psi)(t).
$$
Thus $\Phi(\psi,\psi_0,u)=0$ if and only if $\psi$
is a solution with data $(\psi_0,u)$.
Since $H^2$ is an algebra and $\cN$ is a polynomial in $\psi$
and $\bar\psi$, the map $\cN:H^2\to H^2$ is $C^\infty$ over $\R$.
Moreover, $B_0$ is continuous and bilinear by
Lemma~\ref{L:3.3}. Hence $\Phi$ is also
$C^\infty$ over $\R$.

To apply the implicit function theorem, we show that
$D_\psi\Phi(\psi^0,\psi_0^0,u^0)$ is a bounded
isomorphism of~$X_T$. We have
$$
\bigl(D_\psi\Phi(\psi^0,\psi_0^0,u^0)h\bigr)(t)
=h(t)
+i\int_0^t U(t-r)
\cN'(\psi^0(r))h(r)\,\dd r
+iB_0(u^0,h)(t).
$$
Let $M$ be an upper bound for the operator norms of
$\cN'(\psi^0(r))$ on $H^2$, and choose a finite
partition of~$[0,T]$ into intervals $[a,b]$ such that
$$
M(b-a)+C_T\|u^0\|_{L^2(a,b)}\le\frac12.
$$
For any $f\in X_T$, we solve
$D_\psi\Phi(\psi^0,\psi_0^0,u^0)h=f$
successively on these intervals.
On each interval, the contribution from $[0,a]$
is already known, and the sum of the two integral
operators over $[a,t]$ has norm at most~$1/2$
by Lemma~\ref{L:3.3} and the choice of the partition.
A standard contraction argument therefore gives a unique $h\in X_T$ and $
\|h\|_{X_T}\le C\|f\|_{X_T}$. 
This proves the required~invertibility.

Since $\Phi(\psi^0,\psi_0^0,u^0)=0$, the implicit
function theorem gives a neighborhood $\mathcal V$ of
$(\psi_0^0,u^0)$ and a real~$C^\infty$ map
$\Psi:\mathcal V\to X_T$ with $\Psi(\psi_0^0,u^0)=\psi^0$
such that
$$
\Phi(\Psi(\psi_0,u),\psi_0,u)=0,
\qquad (\psi_0,u)\in\mathcal V.
$$
By Step~2, $\Psi(\psi_0,u)$ is the unique solution
with data $(\psi_0,u)$.
Differentiating this identity, we see that~$v=D\Psi(\psi_0^0,u^0)(h_0,w)$ satisfies
\eqref{E:3.3} in integral form.
Since the right-hand side of the linearized equation
belongs to $L^2(0,T;L^2)$, this identity gives
$v\in W^{1,2}(0,T;L^2)$, and the equation holds
almost everywhere in time as an identity in $L^2$.
Uniqueness follows from the invertibility proved above.
\end{proof}

Clearly, every solution $\psi\in X_T$ conserves the $L^2$ norm.
The equation is also time reversible, that is, if~$\psi$ is a
solution on $[0,T]$ with control $u$, then $t\mapsto\overline{\psi(T-t)}$
is a solution with control $t\mapsto u(T-t)$.

When $\kappa=0$, equation~\eqref{E:0.5} is linear,
and a standard continuation argument gives the
following global well-posedness result.

\begin{corollary}\label{C:3.5}
Assume $\kappa=0$. For any $T>0$, $\psi_0\in H^2$,
and $u\in L^2(0,T;\R^5)$, equation~\eqref{E:0.5}
has a unique solution in $X_T$ with data $(\psi_0,u)$.
\end{corollary}

The smoothing and well-posedness results above also hold
with $Q_c,Q_s$ replaced by $Q_c+R_c,Q_s+R_s$ for arbitrary
$R_c,R_s\in H^2(\T;\R)$.
Indeed, multiplication by $R_c$ and $R_s$ is bounded
on $H^2$, so their Duhamel terms satisfy the estimate
of Lemma~\ref{L:3.3}.

\section{Exact controllability}\label{S:4}

In this section, we assume that $V=0$ and prove
Theorem~\ref{T:B}. We first give more general
conditions on the two additional control potentials.
We then prove local exact controllability near $\phi_0$
and combine it with~Theorem~\ref{T:A}.

\subsection{Admissible control potentials}\label{SS:4.1}

We consider equation~\eqref{E:0.5} with
$Q=(\mathbf1,\cos x,\sin x,P_c,P_s)$.
We assume that $(P_c,P_s)$ is admissible in the following sense.

\begin{definition}\label{D:4.1}
A pair $(P_c,P_s)\in H^1(\T;\R)^2$ with zero means
is called \emph{admissible} if the following conditions hold:
\begin{itemize}
\item[\hypertarget{admissible-a}{\textup{(a)}}]
for $j=c,s$, the map $h\mapsto P_j(h-h(0))$
is a bounded linear operator on $H^2$;
\item[\hypertarget{admissible-b}{\textup{(b)}}] the matrices
$$
M_n=n^2\begin{pmatrix}
\Real\widehat P_c(n)&\Real\widehat P_s(n)\\
\Imag\widehat P_c(n)&\Imag\widehat P_s(n)
\end{pmatrix},\qquad n\ge1,
$$
satisfy
\begin{equation}\label{E:4.1}
\sup_{n\ge2}\big(\|M_n\|+\|M_n^{-1}\|\big)<\infty.
\end{equation}
\end{itemize}
\end{definition}
Condition~\textup{(\hyperlink{admissible-a}{a})} is the multiplication
estimate used in Section~\ref{S:3} (cf.~Lemma~\ref{L:3.1}). We use
condition~\textup{(\hyperlink{admissible-b}{b})} to control the sine
and cosine modes at each frequency $n\ge2$ and to estimate the
controls in $L^2$. Since~$\cos x$ and $\sin x$ are already among the control potentials,
no condition on $M_1$ is needed.

Condition~\textup{(\hyperlink{admissible-b}{b})} also determines
the Sobolev regularity of $P_c$ and $P_s$.
 Indeed,  
the uniform bounds on~$M_n$ and $M_n^{-1}$ imply that
\begin{equation}\label{E:4.2}
c \le n^2 |\widehat P_j(n)|\le C,\qquad n\ge2.
\end{equation}
Thus $P_j\in H^r$ for every $r<3/2$, but
$P_j\notin H^{3/2}$. For the pair $(Q_c,Q_s)$, we have $M_n=\operatorname{diag}(1/2,-1/2)$
for every $n\ge1$, so admissibility follows from Lemma~\ref{L:3.1}.

Since $\mathbf1$ is a control potential, we can absorb the means
of $P_c$ and $P_s$ into the first control. With this normalization, the set of $H^2$ perturbations
that preserve admissibility is open and dense.

\begin{lemma}\label{L:4.2}
Let $P^0=(P_c^0,P_s^0)$ be an admissible pair. For
$R=(R_c,R_s)\in H^2(\T;\R)^2$, let
$$
m_j=\frac1{2\pi}\int_\T R_j(x)\,\dd x,\qquad
P_j=P_j^0+R_j-m_j,\qquad j=c,s.
$$
The set of pairs $R$ for which $(P_c,P_s)$ is admissible
is open and dense in $H^2(\T;\R)^2$ and contains $(0,0)$.
\end{lemma}

The proof is given in Appendix~\ref{SS:B}.

\begin{proposition}\label{P:4.3}
Let $(P_c,P_s)$ be an admissible pair. Then the conclusions of
Lemma~\ref{L:3.3}, Proposition~\ref{P:3.4}, and
Corollary~\ref{C:3.5} hold for
equation~\eqref{E:0.5} with
$Q=(\mathbf1,\cos x,\sin x,P_c,P_s)$.
\end{proposition}
\begin{proof}
It suffices to prove Lemma~\ref{L:3.3} with
$P_c,P_s$ in place of $Q_c,Q_s$.
For $j=c,s$, let us write
$$
P_jh=P_j(h-h(0))+h(0)P_j.
$$
By condition~\textup{(\hyperlink{admissible-a}{a})},
the first term is bounded in $H^2$ by $C\|h\|_2$.
For the second term, the upper bound in~\eqref{E:4.2}, together
with the relation $\widehat P_j(-n)=\overline{\widehat P_j(n)}$ for the
real-valued function $P_j$, gives
$$
\sup_{n\ne0}|n|^2|\widehat P_j(n)|<\infty.
$$
The zero-mean assumption ensures that the zeroth
Fourier coefficient of the Duhamel integral of
$h(0)P_j$ vanishes.
For the remaining coefficients, the bound above
and Lemma~\ref{L:3.2} give the required $H^2$
estimate, as in the proof of Lemma~\ref{L:3.3}.
Applying this estimate on shorter intervals
also implies continuity in $H^2$.

The conclusions of Proposition~\ref{P:3.4} and
Corollary~\ref{C:3.5} then follow
by the same arguments.
\end{proof}

For these control fields, we also use the notation
$\Theta(\psi_0,T)$ for the controls whose solutions
exist in~$C([0,T];H^2)$, and write
$$
\cR_t(\psi_0,u)=\Psi(\psi_0,u)(t).
$$
Uniqueness makes this notation consistent with the
one used for the smooth control fields in Section~\ref{S:1}.

\subsection{Local exact controllability}\label{SS:4.2}

We prove the following local exact controllability
result by establishing controllability of the
linearized equation and then applying the inverse
mapping theorem, as in~\cite[Section~3.3]{BL}.
\begin{theorem}\label{T:4.4}
Assume \eqref{E:0.6}, let $(P_c,P_s)$ be an admissible
pair, and set $Q=(\mathbf1,\cos x,\sin x,P_c,P_s)$.
For~any~$T>0$, there is $\delta>0$ such that, for any $\psi_0,\psi_1\in H^2\cap\bS$ satisfying
$$
\|\psi_j-\phi_0\|_2<\delta,\qquad j=0,1,
$$
there are a control $u\in L^2(0,T;\R^5)$ and
a solution $\psi\in C([0,T];H^2)$ of~\eqref{E:0.5}
such that
$$
\psi(0)=\psi_0,\qquad \psi(T)=\psi_1.
$$
\end{theorem}

Let us fix an admissible pair and linearize the equation
at the constant state $\phi_0$. We put
$$
\kappa_0=\kappa(2\pi)^{-p},\qquad u^0=(-\kappa_0,0,0,0,0).
$$
Since $\partial_x^2\phi_0=0$ and
$\cN(\phi_0)=\kappa_0\phi_0$, the constant function
$\psi(t,x)\equiv\phi_0$ is a solution of
equation~\eqref{E:0.5} with control $u^0$. 
The tangent space to $H^2\cap\bS$ at $\phi_0$ and
the projection onto it are given by
$$
E=\{v\in H^2:\Real\hat v(0)=0\},\qquad \Pi\psi=\psi-\Real\hat\psi(0).
$$
We shall first prescribe the projection of the
state at time $T$ onto $E$. Near $\phi_0$, the remaining
real constant part is positive and can be
recovered from the $L^2$ norm.
Indeed, orthogonality gives
\begin{equation}\label{E:4.3}
\|\psi\|_{L^2}^2=\|\Pi\psi\|_{L^2}^2+2\pi\big(\Real\hat\psi(0)\big)^2.
\end{equation}

Recalling that $\beta=2p\kappa_0$, we have
$$
\cN'(\phi_0)v+\lan u^0,Q\ran v
=\kappa_0(p+1)v+\kappa_0p\bar v-\kappa_0v
=\beta\Real v.
$$
Consequently, the equation linearized at $(\phi_0,u^0)$ is
\begin{equation}\label{E:4.4}
i\partial_t v=-\partial_x^2v+\beta\Real v
+\phi_0\lan w(t),Q(x)\ran,\qquad v(0)=0.
\end{equation}
We use the two additional potentials to control
the modes $n\ge2$ by a moment problem.
The condition $\beta>-4$ ensures that $n^2+\beta$
is positive for all these modes.
The first mode is controlled directly by
$\cos x$ and $\sin x$, so no sign condition
on $1+\beta$ is needed.
The constant potential controls the zero-mode phase.

The following proposition gives a bounded right
inverse for the map $w\mapsto v(T)$ and thus
proves exact controllability of the linearized
equation.

\begin{proposition}\label{P:4.5}
Assume \eqref{E:0.6}, let $(P_c,P_s)$ be admissible, and let $T>0$.
\begin{itemize}
\item[(i)] For any $w\in L^2(0,T;\R^5)$, equation~\eqref{E:4.4} has a unique solution $v\in X_T$, and $v(t)\in E, t\in[0,T]$.
\item[(ii)] There is a bounded real-linear operator $\cM:E\to L^2(0,T;\R^5)$ such that the solution corresponding to~$w=\cM v^1$ satisfies $v(T)=v^1$ for every $v^1\in E$.
\end{itemize}
\end{proposition}
\begin{proof}
\emph{Step~1. The Fourier equations.}
Let $w\in L^2(0,T;\R^5)$. By
Proposition~\ref{P:4.3}, the function
$v=D\Psi(\phi_0,u^0)(0,w)\in X_T$ is the unique
solution of~\eqref{E:4.4}.
We denote the Fourier coefficients of $\Real v$ and~$\Imag v$ by $a_n$ and $b_n$, respectively, and set
$$
q_n=\sum_{j=1}^5w_j\hat Q_j(n).
$$
Thus $\hat v(n)=a_n+ib_n$, with
$a_{-n}=\overline{a_n}$ and $b_{-n}=\overline{b_n}$.
Taking Fourier coefficients of the real and imaginary
parts of~\eqref{E:4.4}, we obtain
\begin{equation}\label{E:4.5}
\dot a_n=n^2b_n,\qquad
\dot b_n=-(n^2+\beta)a_n-\phi_0q_n
\end{equation}
almost everywhere on $[0,T]$, with $a_n(0)=b_n(0)=0$.
Since all potentials except $\mathbf1$ have zero mean,
$q_0=w_1$, and hence
$$
a_0(t)=0,\qquad
b_0(t)=-\phi_0\int_0^t w_1(r)\,\dd r.
$$
This proves (i). For $n\ge2$, let us set
$$
\omega_n=n\sqrt{n^2+\beta},\qquad
\lambda_n=\frac{\omega_n}{n^2}
=\sqrt{1+\beta/n^2}.
$$
Since $\beta>-4$, the numbers $\lambda_n$ are positive
and bounded above and away from zero.
The functions $z_n^\pm=\lambda_na_n\pm ib_n$ satisfy
$$
\dot z_n^\pm=\mp i\omega_nz_n^\pm\mp i\phi_0q_n,
\qquad z_n^\pm(0)=0.
$$
Solving the diagonalized equations, we obtain
(cf.~\cite[Section~3.3]{BL})
\begin{equation}\label{E:4.6}
z_n^\pm(T)=\mp i\phi_0e^{\mp i\omega_nT}
\int_0^T e^{\pm i\omega_nt}q_n(t)\,\dd t.
\end{equation}

\emph{Step~2. The moment problem.}
For $F(x)=x\sqrt{x^2+\beta}$, we~have
$$
F'(x)=\sqrt{x^2+\beta}
+\frac{x^2}{\sqrt{x^2+\beta}}\ge2x,\qquad x\ge2.
$$
Integrating over $[n,n+1]$ for $n\ge2$, we obtain
\begin{equation}\label{E:4.7}
\omega_{n+1}-\omega_n\ge2n+1\longrightarrow\infty.
\end{equation}
Moreover, $\omega_2=2\sqrt{4+\beta}>0$.
The divergence of the frequency gaps allows us
to solve the moment problem on every interval
of positive length.
Applying \cite[Appendix~B, Corollary~1]{BL} to the
sequence $0,\omega_2,\omega_3,\ldots$, with zero
moment prescribed at frequency zero, we obtain a bounded
real-linear operator
$$
\cJ:\ell^2(\{2,3,\ldots\};\C)\longrightarrow L^2(0,T;\R),
$$
such that $w=\cJ d$ satisfies
$$
\int_0^T w(t)e^{i\omega_nt}\,\dd t=d_n,\qquad n\ge2,
\qquad
\|w\|_{L^2(0,T)}\le C_{T,\beta}\|d\|_{\ell^2}.
$$

It remains to choose the right-hand sides $d_n$
of these moment equations so that the corresponding
controls steer the solution to the prescribed target.

\emph{Step~3. Choice of the controls for $n\ge2$.}
Let $v^1\in E$ be the prescribed final state, and
denote by~$a_n^1,b_n^1$ the Fourier coefficients of
its real and imaginary parts. For $n\ge2$, set
$z_n^{\pm,1}=\lambda_na_n^1\pm ib_n^1$.
Since the sequence $(\lambda_n)$ is bounded,
the $H^2$ regularity of~$v^1$ gives
\begin{equation}\label{E:4.8}
\sum_{n=2}^\infty n^4
\bigl(|z_n^{+,1}|^2+|z_n^{-,1}|^2\bigr)
\le C\|v^1\|_2^2.
\end{equation}
The factor $n^2$ appearing in the definition of
$M_n$ is chosen precisely to match these
$H^2$ weights.

We write
$$
w_c=w_4,\qquad w_s=w_5.
$$
Since the first three control potentials have zero
$n$th Fourier coefficients for $n\ge2$, we have
$$
q_n(t)
=\widehat P_c(n)w_c(t)+\widehat P_s(n)w_s(t),
\qquad n\ge2.
$$
In view of~\eqref{E:4.6}, to obtain
$z_n^\pm(T)=z_n^{\pm,1}$ we need to prescribe the
moments of $w_c$ and $w_s$ at
the frequencies~$\omega_n$.
Let
$$
d_n^\pm=\pm\frac{in^2}{\phi_0}
e^{\pm i\omega_nT}z_n^{\pm,1}
$$
and  
$$
d_{j,n}
:=\int_0^T e^{i\omega_nt}w_j(t)\,\dd t,
\qquad j=c,s.
$$
Since the controls are real-valued, their moments at
$-\omega_n$ are $\overline{d_{j,n}}$.
Substituting these moments into
\eqref{E:4.6} and conjugating the
equation for $z_n^-(T)$, we find that
$z_n^\pm(T)=z_n^{\pm,1}$ is equivalent to
\begin{equation}\label{E:4.9}
\begin{pmatrix}
n^2\widehat P_c(n)&
n^2\widehat P_s(n)\\
n^2\overline{\widehat P_c(n)}&
n^2\overline{\widehat P_s(n)}
\end{pmatrix}
\begin{pmatrix}d_{c,n}\\d_{s,n}\end{pmatrix}
=
\begin{pmatrix}d_n^+\\\overline{d_n^-}\end{pmatrix}.
\end{equation}
The matrix on the left is
$$
\begin{pmatrix}1&i\\1&-i\end{pmatrix}M_n.
$$
By condition~\textup{(\hyperlink{admissible-b}{b})},
its inverse is uniformly bounded for $n\ge2$.
Solving~\eqref{E:4.9}, we obtain
$$
|d_{c,n}|^2+|d_{s,n}|^2
\le C\bigl(|d_n^+|^2+|d_n^-|^2\bigr).
$$
Using~\eqref{E:4.8}, we obtain
$$
\sum_{n=2}^\infty
\bigl(|d_{c,n}|^2+|d_{s,n}|^2\bigr)
\le C\|v^1\|_2^2.
$$
Thus $d_c,d_s\in\ell^2$, and the moment operator
$\cJ$ from Step~2 gives real controls
$$
w_4=\cJ d_c,\qquad w_5=\cJ d_s
$$
having these prescribed moments. Since $\lambda_n>0$,
we obtain $a_n(T)=a_n^1$ and $b_n(T)=b_n^1$ for
every $n\ge2$.

\emph{Step~4. The first and zero Fourier modes.}
We now choose the controls $w_2,w_3$ to obtain
the prescribed first Fourier mode. Let us set
$r=t/T$ and
$$
A(t)=a_1^1(3r^2-2r^3)+Tb_1^1(r^3-r^2).
$$
Then
$$
A(0)=A'(0)=0,\qquad A(T)=a_1^1,\qquad A'(T)=b_1^1.
$$
We define
$$
R(t)=-\frac{A''(t)+(1+\beta)A(t)}{\phi_0}
-\widehat P_c(1)w_4(t)-\widehat P_s(1)w_5(t),
$$
and choose
$$
w_2(t)=2\Real R(t),\qquad
w_3(t)=-2\Imag R(t).
$$
Since $\widehat{\cos x}(1)=1/2$ and
$\widehat{\sin x}(1)=-i/2$, we have
$$
q_1(t)=\frac{w_2(t)-iw_3(t)}2
+\widehat P_c(1)w_4(t)+\widehat P_s(1)w_5(t)
=-\frac{A''(t)+(1+\beta)A(t)}{\phi_0}.
$$
Thus $a_1=A$ and $b_1=A'$ solve
\eqref{E:4.5} for $n=1$, and
$$
a_1(T)=a_1^1,\qquad b_1(T)=b_1^1.
$$
The controls $w_2,w_3$ do not affect the modes
$n\ge2$.

It remains to control the zero Fourier mode.
Since $v^1\in E$, we have $a_0^1=0$, while
$b_0^1=\Imag\widehat v^1(0)$.
We~therefore set
$$
w_1(t)=-\frac{b_0^1}{\phi_0T}.
$$
This gives
$$
a_0(T)=a_0^1=0,\qquad b_0(T)=b_0^1,
$$
and hence all nonnegative Fourier coefficients of
$\Real v(T)$ and $\Imag v(T)$ agree with those of
the target. By conjugate symmetry, $v(T)=v^1$.

Finally, the construction is real-linear in $v^1$,
and the estimates from Step~3 and the definition of $A$ give
$$
\|w\|_{L^2(0,T;\R^5)}
\le C_{T,\beta,P_c,P_s}\|v^1\|_2.
$$
Thus $\cM v^1=w$ is the required bounded right inverse.
\end{proof}

\begin{proof}[Proof of Theorem~\ref{T:4.4}]
Let us fix $T>0$ and denote by $\cM$ the operator
given by Proposition~\ref{P:4.5}.
By Proposition~\ref{P:4.3}, the map
$$
F_T(\psi_0,v)
=\Pi\big[\Psi(\psi_0,u^0+\cM v)(T)\big]
$$
is well defined and $C^\infty$ over $\R$ in a
neighborhood of $(\phi_0,0)$ in $H^2\times E$.
Moreover, Proposition~\ref{P:4.5} implies
$$
F_T(\phi_0,0)=0,\qquad
D_vF_T(\phi_0,0)=\Id_E.
$$
Consider now the map
$$
\mathcal F_T(\psi_0,v)
=\bigl(\psi_0,F_T(\psi_0,v)\bigr)
$$
from a neighborhood of $(\phi_0,0)$ in $H^2\times E$
to $H^2\times E$.
Its derivative at $(\phi_0,0)$ is triangular with
identity diagonal blocks and is therefore a bounded
isomorphism.
The inverse mapping theorem gives a map~$(\psi_0,y)\mapsto v(\psi_0,y)\in E$, of class
$C^\infty$ over $\R$, defined near $(\phi_0,0)$
and satisfying
$$
v(\phi_0,0)=0,\qquad
F_T(\psi_0,v(\psi_0,y))=y.
$$

For $\psi_0,\psi_1\in H^2\cap\bS$ sufficiently close
to $\phi_0$, let us take
$$
v=v(\psi_0,\Pi\psi_1),\qquad
u=u^0+\cM v,\qquad
\psi=\Psi(\psi_0,u).
$$
By construction, $\Pi\psi(T)=\Pi\psi_1$.
It remains to prove that the real constant parts
also agree.
Conservation of the $L^2$ norm gives
$\|\psi(T)\|_{L^2}=\|\psi_1\|_{L^2}=1$.
Together with~\eqref{E:4.3}, this implies
$$
\bigl(\Real\widehat{\psi(T)}(0)\bigr)^2
=\bigl(\Real\widehat{\psi_1}(0)\bigr)^2.
$$
By continuity of the maps constructed above,
$\psi(T)\to\phi_0$ in $H^2$ as
$(\psi_0,\psi_1)\to(\phi_0,\phi_0)$.
Since $\widehat\phi_0(0)=(2\pi)^{-1/2}>0$,
both real parts are positive for $\psi_0,\psi_1$
sufficiently close to $\phi_0$.
They are therefore equal, and~$\psi(T)=\psi_1$.
Choosing $\delta>0$ sufficiently small, we complete
the proof.
\end{proof}

\subsection{Global exact controllability}\label{SS:4.3}

We use the three smooth potentials to bring
the initial state and the complex conjugate
of the target close to $\phi_0$ in $H^2$.
The two additional potentials are used only
for the local exact controllability.
We then complete the construction by time
reversibility. Since the approximate and local
exact results hold at every given positive
time, we obtain the following global controllability result.

\begin{theorem}\label{T:4.6}
Assume \eqref{E:0.6} and let $(P_c,P_s)$ be an admissible pair.
For any $T>0$ and any initial and target states
$\psi_0,\psi_1\in H^2\cap\bS$, there are a control
$u\in L^2(0,T;\R^5)$ and a solution
$\psi\in C([0,T];H^2\cap\bS)$ of~\eqref{E:0.5}
such that
$$
\psi(0)=\psi_0,\qquad \psi(T)=\psi_1.
$$
\end{theorem}
\begin{proof}
Let $\tau=T/3$ and choose $\delta>0$ as in
Theorem~\ref{T:4.4} for time $\tau$.
By Theorem~\ref{T:A}, applied with $s=2$ and
the field $(\mathbf1,\cos x,\sin x)$, there are controls
$u_1,v\in L^2(0,\tau;\R^5)$, with their last two
components equal to zero, such that the states
$$
h_0=\cR_\tau(\psi_0,u_1),\qquad
h_1=\overline{\cR_\tau(\overline{\psi_1},v)}
$$
satisfy $\|h_j-\phi_0\|_2<\delta$, $j=0,1$.
The estimate for $h_1$ uses the fact that $\phi_0$ is real.

By Theorem~\ref{T:4.4}, there is a control
$u_2\in L^2(0,\tau;\R^5)$ such that
$\cR_\tau(h_0,u_2)=h_1$.
By time reversibility, the control
$u_3(r)=v(\tau-r)$ satisfies
$\cR_\tau(h_1,u_3)=\psi_1$.
Concatenating $u_1,u_2,u_3$ completes the proof.
\end{proof}

\begin{proof}[Proof of Theorem~\ref{T:B}]
Let $\cO_0$ be the set of perturbations given by
Lemma~\ref{L:4.2} for $P^0=(Q_c,Q_s)$, and set
$$
\cO=(Q_c,Q_s)+\cO_0\subset\mathcal Q.
$$
This set is open and dense in $\mathcal Q$ and contains
$(Q_c,Q_s)$.
For $(P_c,P_s)\in\cO$, let us write
$P_j=Q_j+R_j$, $j=c,s$, and denote the means of
$R_c,R_s$ by $m_c,m_s$.
The pair $(P_c-m_c,P_s-m_s)$ is admissible.
Let~$v=(v_1,\dots,v_5)\in L^2(0,T;\R^5)$ be a control
given by Theorem~\ref{T:4.6} for this pair, and set
$$
u=(v_1-m_cv_4-m_sv_5,\,v_2,\,v_3,\,v_4,\,v_5).
$$
The same solution satisfies the equation with control
$u$ and the pair $(P_c,P_s)$.
This proves the theorem.
\end{proof}

\begin{remark}\label{R:4.7}
Let us explain how resonances arise in the linearization
at $(\phi_0,u^0)$. By~\eqref{E:4.5}, when $w=0$,
the coefficient $a_n$ satisfies
$$
a_n''=-n^2(n^2+\beta)a_n.
$$
Under condition~\eqref{E:0.6}, all modes $n\ge2$ are
oscillatory, with distinct positive frequencies. If
$n^2+\beta<0$, the $n$th mode instead has one growing
and one decaying exponential, with exponents
$\pm n\sqrt{-n^2-\beta}$.

For two distinct modes $2\le n<m$, the exponents coincide
when $\beta=-(n^2+m^2)$, which lies outside the range
of~\eqref{E:0.6}, and their common values are $\pm nm$.
At each shared exponent, the two controls multiplying~$P_c$ and $P_s$ provide only two real moments, whereas
the sine and cosine components of the two modes require
four. The other three controls do not act on these modes.
Thus the linearization at $(\phi_0,u^0)$ is not exactly
controllable at these values of $\beta$.
This does not rule out controllability of the nonlinear~equation.
\end{remark}

\appendix
\setcounter{tocdepth}{1}
\addtocontents{toc}{\protect\setcounter{tocdepth}{1}}
\refstepcounter{section}
\section*{Appendix}
\addcontentsline{toc}{section}{Appendix}
\renewcommand{\thesubsection}{\Alph{subsection}}
\setcounter{subsection}{0}

In this section, we prove two auxiliary lemmas used in the
previous sections.

\subsection{Choice of a periodic phase}\label{SS:A}

For each odd integer $N\ge1$, let $E_N$ denote
the $N$-dimensional space of real trigonometric
polynomials of degree at most $(N-1)/2$.
The following lemma gives the choice of phase
used in Lemma~\ref{L:2.4}.

\begin{lemma}\label{L:A.1}
Let $N\ge3$ be an odd integer and $f_0,\dots,f_{N-1}\in H^1(\T;\C)$.
Assume that, for every $x\in\T$, at least three of the values
$f_0(x),\dots,f_{N-1}(x)$ are nonzero. Then, for Lebesgue almost
every $\theta\in E_N$,
$$
\sum_{k=0}^{N-1}f_k(x)e^{i\theta(x+2\pi k/N)}\ne0
\qquad\text{for every }x\in\T.
$$
\end{lemma}

\begin{proof}
\emph{Step~1. An estimate for independent angles.}
Let us first allow the $N$ angles in the sum to vary independently.
For $a_0,\dots,a_{N-1}\in\C$, we shall prove that
\begin{equation}\label{E:A.1}
\left|\left\{\alpha\in[0,2\pi]^N:
\left|\sum_{k=0}^{N-1}a_ke^{i\alpha_k}\right|\le r
\right\}\right|\le C_c r^{3/2},\qquad r>0,
\end{equation}
provided that at least three  coefficients $a_k$ have modulus
at least $c>0$. Here the outer bars denote Lebesgue measure,
and $C_c$ depends only on $c$ and $N$.
This estimate can be derived from the
Fourier estimates in~\cite[proof of Lemma~3.2]{CM12}
or~\cite[Example~2.6.9]{Gra14}. We give a detailed proof for~completeness.

We identify $\C$ with $\R^2$ and set
$$
F(\alpha)=\sum_{k=0}^{N-1}a_ke^{i\alpha_k},\qquad
\Phi(\xi)=\int_{[0,2\pi]^N}e^{i\xi\cdot F(\alpha)}\,\dd\alpha,
\qquad \xi\in\R^2.
$$
To estimate $\Phi$, let us consider the function
$$
\chi(t)=\frac1{2\pi}\int_0^{2\pi}e^{it\cos\alpha}\,\dd\alpha
=\frac1\pi\int_0^\pi e^{it\cos\alpha}\,\dd\alpha,
\qquad t\ge0.
$$
We shall show that
\begin{equation}\label{E:A.2}
|\chi(t)|\le C(1+t)^{-1/2},\qquad t\ge0.
\end{equation}
Since $|\chi(t)|\le1$, it suffices to consider $t\ge1$.
Let $\eta=t^{-1/2}$. The integrals over $[0,\eta]$ and
$[\pi-\eta,\pi]$ are each bounded in modulus by $\eta$.
On $[\eta,\pi-\eta]$, integration by parts gives
$$
\left|\int_\eta^{\pi-\eta}e^{it\cos\alpha}\,\dd\alpha\right|
\le \frac{2}{t\sin\eta}
+\frac1t\int_\eta^{\pi-\eta}
\frac{|\cos\alpha|}{\sin^2\alpha}\,\dd\alpha
\le \frac{4}{t\sin\eta}
\le \frac{C}{t\eta}.
$$
Combining the estimates on the three subintervals of $[0,\pi]$,
we obtain
$$
|\chi(t)|\le \frac1\pi\left(2\eta+\frac{C}{t\eta}\right)
\le Ct^{-1/2},\qquad t\ge1,
$$
which proves~\eqref{E:A.2}. Let us now write $a_k=|a_k|e^{i\gamma_k}$ and
$\xi=|\xi|e^{i\gamma}$. Since
$$
\xi\cdot(a_ke^{i\alpha})
=\Real\bigl(\overline{\xi}\,a_ke^{i\alpha}\bigr)
=|a_k||\xi|\cos(\alpha+\gamma_k-\gamma),
$$
periodicity gives
$$
\int_0^{2\pi}e^{i\xi\cdot(a_ke^{i\alpha})}\,\dd\alpha
=2\pi\chi(|a_k||\xi|).
$$
It follows from Fubini's theorem that
\begin{equation}\label{E:A.3}
\Phi(\xi)=(2\pi)^N\prod_{k=0}^{N-1}\chi(|a_k||\xi|).
\end{equation}
Applying~\eqref{E:A.2} to three factors for which
$|a_k|\ge c$ and using $|\chi|\le1$ for the others, we obtain
\begin{equation}\label{E:A.4}
|\Phi(\xi)|\le C_N(1+c|\xi|)^{-3/2},\qquad \xi\in\R^2.
\end{equation}

To prove~\eqref{E:A.1}, let us choose a nonnegative
function $\zeta\in C_c^\infty(\R^2)$ equal to $1$ on the unit disk.
Then~$\zeta(F(\alpha)/r)=1$ whenever $|F(\alpha)|\le r$.
By Fourier inversion and Fubini's theorem,
$$
\bigl|\{\alpha\in[0,2\pi]^N:|F(\alpha)|\le r\}\bigr|
\le \int_{[0,2\pi]^N}\zeta\bigl(F(\alpha)/r\bigr)\,\dd\alpha =\frac{r^2}{(2\pi)^2}
\int_{\R^2}\widehat\zeta(r\xi)\Phi(\xi)\,\dd\xi,
$$
where $\widehat\zeta(\xi)=\int_{\R^2}e^{-i\xi\cdot z}\zeta(z)\,\dd z$.
By~\eqref{E:A.4} and the change of variables $\eta=r\xi$,
\begin{equation}\label{E:A.5}
r^2\int_{\R^2}|\widehat\zeta(r\xi)\Phi(\xi)|\,\dd\xi
\le C_N(r/c)^{3/2}
\int_{\R^2}|\widehat\zeta(\eta)|\,|\eta|^{-3/2}\,\dd\eta.
\end{equation}
The last integral is finite, since $|\eta|^{-3/2}$ is integrable
near the origin and $\widehat\zeta$ is rapidly decreasing.
This proves~\eqref{E:A.1}.

\emph{Step~2. An estimate in the phase space.}
Let $\tau=2\pi/N$ and set
$$
\cZ(x,\theta)=\sum_{k=0}^{N-1}f_k(x)e^{i\theta(x+k\tau)}.
$$
We endow $E_N$ with the $L^2$ norm and the corresponding
Lebesgue measure. For $R>0$, let $B_R\subset E_N$ be the
closed ball of radius $R$ centered at the origin. In this step,
we prove that
\begin{equation}\label{E:A.6}
\bigl|\{\theta\in B_R:|\cZ(x,\theta)|\le r\}\bigr|
\le C_Rr^{3/2},\qquad x\in\T,\quad r>0,
\end{equation}
where $C_R>0$ is independent of $x$ and $r$.

At each point of $\T$, three of the functions $f_k$ are nonzero.
Their moduli have a positive lower bound on a neighborhood of
that point. By taking a finite cover of $\T$, we obtain $c>0$
such that, at every $x\in\T$, at least three of the values
$|f_k(x)|$ are at least $c$.

For $x\in\T$, let
$$
A_x:E_N\longrightarrow\R^N,\qquad
A_x\theta=\bigl(\theta(x+k\tau)\bigr)_{k=0}^{N-1}.
$$
By trigonometric interpolation, $A_x$ is an isomorphism.
Its matrix depends continuously on $x$, so compactness of $\T$
gives $|\det A_x|^{-1}\le C$ for every $x\in\T$.

Let us set
$$
D_{x,r}=\left\{\alpha\in A_x(B_R):
\left|\sum_{k=0}^{N-1}f_k(x)e^{i\alpha_k}\right|\le r\right\}.
$$
The change of variables $\alpha=A_x\theta$ gives
\begin{equation}\label{E:A.7}
\bigl|\{\theta\in B_R:|\cZ(x,\theta)|\le r\}\bigr|
=|\det A_x|^{-1}|D_{x,r}|\le C|D_{x,r}|.
\end{equation}

Let $M_R=\sup_{\theta\in B_R}\|\theta\|_{L^\infty}$, which is
finite since $E_N$ is finite dimensional. Each coordinate of
$A_x\theta$ is a value of $\theta$, so
$$
D_{x,r}\subset A_x(B_R)\subset[-M_R,M_R]^N
\qquad\text{for every }x\in\T.
$$
Let us cover $[-M_R,M_R]^N$ by $L_R$ cubes of the form
$2\pi\ell+[0,2\pi]^N$, $\ell\in\Z^N$.
By periodicity in the angles, we may apply~\eqref{E:A.1}
with $a_k=f_k(x)$ on each cube $Q$ of the covering. Thus
\begin{equation}\label{E:A.8}
|D_{x,r}|\le \sum_Q |D_{x,r}\cap Q|\le L_R C_c r^{3/2}.
\end{equation}
In view of~\eqref{E:A.7}, this proves
\eqref{E:A.6} with $C_R=C L_R C_c$.

\emph{Step~3. Excluding zeros on the whole circle.}
Let us fix $R>0$ and let $\mathcal E_R$ be the set of
$\theta\in B_R$ for which the function
$\cZ(\cdot,\theta):\T\to\C$ has a zero.
This set is closed by continuity of $\cZ$ and compactness of $\T$.
Let us prove that $|\mathcal E_R|=0$.

Let $h=2\pi/m$ and divide $[0,2\pi]$ into $m$ intervals
$I_j$ of length $h$, with left endpoints $x_j$.
  As $\|\theta'\|_{L^\infty}$ is bounded on $B_R$,
the Cauchy--Schwarz inequality gives
$$
|\cZ(x,\theta)-\cZ(x_j,\theta)|^2
\le K_Rh\left(
h+\sum_{k=0}^{N-1}\int_{I_j}|f_k'(y)|^2\,\dd y
\right)=:r_j^2
$$
for every $x\in I_j$ and $\theta\in B_R$, where $K_R$
is independent of $m$ and $j$.

If $\cZ(\cdot,\theta)$ vanishes in $I_j$, then
$|\cZ(x_j,\theta)|\le r_j$.
Every $\theta\in\mathcal E_R$ satisfies this inequality
for at least one $j$. Hence
$$
\mathcal E_R\subset
\bigcup_{j=0}^{m-1}
\{\theta\in B_R:|\cZ(x_j,\theta)|\le r_j\}.
$$
By subadditivity of Lebesgue measure, \eqref{E:A.6},
and H\"older's inequality,
$$
|\mathcal E_R|
\le C_R\sum_jr_j^{3/2}
\le C_R m^{1/4}\Big(\sum_j r_j^2\Big)^{3/4}
\le C h^{1/2}\longrightarrow0
$$
as $m\to\infty$, since $mh=2\pi$ and
$$
\sum_j r_j^2=K_Rh(2\pi+\sum_{k=0}^{N-1}\|f_k'\|_{L^2}^2).
$$
Hence $|\mathcal E_R|=0$.
Taking the union over positive integers $R$ proves the lemma.
\end{proof}

With only two nonzero terms, the conclusion of
Lemma~\ref{L:A.1} may fail. Indeed, let $f_0=1$,
$f_1=-1$, and $f_k=0$ for all other indices.
For any continuous real periodic phase $\theta$, the function
$x\mapsto\theta(x+2\pi/N)-\theta(x)$ has zero integral over $\T$.
It therefore vanishes at some point, where the two terms cancel.

\subsection{Proof of Lemma~\ref{L:4.2}}\label{SS:B}

Let $\cO_0$ denote the set of perturbations $R=(R_c,R_s)$
in Lemma~\ref{L:4.2} for which the resulting pair
$P=(P_c,P_s)$ is admissible. We write
$\|R\|=\|R_c\|_2+\|R_s\|_2$.
Since $P^0$ is admissible, $(0,0)\in\cO_0$.
Adding an $H^2$ function to either potential preserves
condition~\textup{(\hyperlink{admissible-a}{a})},
because $H^2$ is an algebra.
We therefore only need to consider
condition~\textup{(\hyperlink{admissible-b}{b})}.

Subtracting the means does not change the nonzero
Fourier coefficients. Thus, with the notation of
Definition~\ref{D:4.1},
$$
M_n(P)=M_n(P^0)+M_n(R),\qquad n\ge2.
$$
Since $R_c,R_s\in H^2$, the sequences
$(n^2\widehat R_j(n))_{n\ge2}$ belong to $\ell^2$.
In particular,
$$
M_n(R)\longrightarrow0,\qquad
\sup_{n\ge2}\|M_n(R)\|\le C\|R\|.
$$
For each fixed $R$, the matrices $M_n(P)$ are therefore
uniformly bounded in $n$. It remains to consider their inverses.
Let us first prove openness. Fix $R\in\cO_0$ and put
$$
K=\sup_{n\ge2}\|M_n(P)^{-1}\|.
$$
If $\|\widetilde R-R\|<(2CK)^{-1}$, then
$$
\sup_{n\ge2}
\|M_n(P)^{-1}M_n(\widetilde R-R)\|
\le CK\|\widetilde R-R\|<\frac12.
$$
The Neumann series shows that the matrices corresponding
to $\widetilde R$ are invertible, with inverses bounded
by~$2K$. Hence $\widetilde R\in\cO_0$, which proves openness.

To prove density, it suffices to approximate any pair $R=(R_c,R_s)$
of real trigonometric polynomials, since such pairs are dense in
$H^2(\T;\R)^2$. Let $N\ge2$ be an upper bound for their degrees
and consider the perturbation $tR$, with $t\in\R$. For $n>N$,
the corresponding matrix is $M_n(P^0)$, so its inverse is
uniformly bounded. For $2\le n\le N$, the determinant
$$
\det\bigl(M_n(P^0)+tM_n(R)\bigr)
$$
is a polynomial in $t$ which is nonzero at $t=0$. Thus, outside
a finite set of values of $t$, all these matrices are invertible
and $tR\in\cO_0$. Taking $t\to1$ outside this set, we obtain
the required approximation.

\section*{Declaration of generative AI use}

The proofs are the author’s. ChatGPT has been used to check the proofs,  to proofread and polish the text,   and to assist in preparing the figures. The author takes full responsibility for the paper.

\section*{Acknowledgments}
This work was partially supported by the National
Natural Science Foundation of China
(Grant No.~12571156) and the NYU Shanghai Boost Fund.

\end{document}